\documentclass[letterpaper, 10pt]{IEEEtran}
\usepackage{amsmath,graphicx,epsfig,color,amsfonts, subfigure,algorithm, psfrag,parskip}
\usepackage{balance}
\usepackage{version,xspace}
\usepackage{framed}
\usepackage[noend]{algorithmic}
\usepackage[sort&compress,square,comma,numbers]{natbib}
\graphicspath{{./fig/}}

\definecolor{shadecolor}{rgb}{0.8,0.8,0.8}

\def\prob{\mathbb{P}}

\def\expt{\mathbb{E}}
\def\real{\mathbb{R}}

\def\naturals{\mathbb{N}}

\def\t{\boldsymbol{t}}

\def\q{\boldsymbol{q}}

\def\t{\boldsymbol{t}}

\newcommand{\until}[1]{\{1,\dots, #1\}}

\newcommand{\subscr}[2]{#1_{\textup{#2}}}
\newcommand{\supscr}[2]{#1^{\textup{#2}}}
\newcommand{\setdef}[2]{\{#1 \; | \; #2\}}
\newcommand{\seqdef}[2]{\{#1\}_{#2}}

\newcommand{\map}[3]{#1: #2 \rightarrow #3}

\newcommand{\subject}{\text{subject to}}
\newcommand{\maximize}{\text{maximize}}
\newcommand{\minimize}{\text{minimize}}

\newcommand{\tran}{^{\top}}

\newcommand\oprocendsymbol{\hbox{$\square$}}
\newcommand\oprocend{\relax\ifmmode\else\unskip\hfill\fi\oprocendsymbol}

\def\tinf{t^{\text{inf}}}

\def\fdi{f^{\dag}}

\def \etal{{\it et al.}}
\def \mc {\mathcal}
\def\fav{\bar{f}}
\def \bs {\boldsymbol}
\title{%Adaptive Information Management Policies \\for Optimal Decision Making in Human-Robot Systems 
Mixed Human-Robot Team Surveillance\\
\Large \it  Integrating Cognitive Modeling with Engineering Design}
\author{Vaibhav~Srivastava$^1$, Amit~Surana$^2$, Miguel~P.~Eckstein$^3$, and Francesco~Bullo$^4$\\ %[-2ex]
{$^1$Department of Mechanical and Aerospace Engineering, Princeton University}\\%[-2ex]
{$^2$United Technologies Research Center, Hartford} \\
{$^3$Department of Psychology and Brain Sciences, UCSB}\\
{$^4$Center for Control, Dynamical Systems and Computation, UCSB}
}

\begin{document}
\maketitle
%\CSMsetup

%\begin{abstract}
%We study the mixed human-robot team design in a system theoretic setting using the context of a surveillance mission. 
%The three key coupled components of a mixed team design are (i) policies for the human operator, (ii) policies to account for erroneous human decisions, and (iii) policies to control the automaton. 
%In this paper, we survey elements of human decision-making, including evidence aggregation, situational awareness, fatigue, and memory effects. We bring together the models for these elements in human decision-making to develop a single coherent model for human decision-making in a two-alternative choice task. We utilize the developed model to design efficient attention allocation policies for the human operator. We propose an anomaly detection algorithm that utilizes potentially erroneous decision by the operator to ascertain an anomalous region among the set of regions surveilled. Finally, we propose a stochastic vehicle routing policy that surveils an anomalous region with high probability. Our mixed team design relies on the certainty-equivalent receding-horizon control framework.
%\end{abstract}

The emergence of mobile and fixed sensor networks operating at different
modalities, mobility, and coverage has enabled access to an unprecedented
amount of information. In a variety of complex and information rich systems,
this information is processed by a human operator~\cite{WMB:09,CD:10}. 
%However, the human operators unable to extract actionable
%%information in a timely manner. 
%Consequently, cooperative
%control and coordination strategies in sensor management have spawned significant research~\cite{CAR-CYS-AT:07, FB-JC-BP:09f}. The human operators often bring knowledge and experience to bear
%in the complex problems where a wrong decision can have lethal consequences. 
%Accordingly,
%the mixed human-robot teams are becoming increasingly important in 
%and have already been deployed in several civilian  as well as military facilities.
%, e.g., surveillance and intelligence gathering.
The inherent inability of humans to handle the plethora of available information
has detrimental effects on their performance and may lead to dire consequences~\cite{TS-MR:11}. 
To alleviate this loss in performance of the human operator, 
the recent National Robotic Initiative~\cite{EG:11} emphasizes collaboration of humans with robotic partners, and envisions 
a symbiotic co-robot that facilitates an efficient interaction of the human operator with the automaton. 
%Given the complex interaction that can arise between  the operator and the automaton, 
An efficient co-robotic partner will enable a better interaction between the automaton and the operator by  exploiting  the operator's strengths while taking into account  their 
inefficiencies, such as erroneous decisions, fatigue and loss of situational awareness.
The design of such co-robots requires algorithms that enable the co-robot to aid the human partner to focus their attention to the pertinent information and direct the automaton to efficiently collect the information. 

In this paper, we focus on the design of mixed human-robot teams. The purpose of mixed teams is to exploit human cognitive abilities in complex missions, and therefore,
%It has been evident that the information overload in these complex missions has had a detrimental effect on the human performance leading to unpleasant consequences~\cite{TS-MR:11}. Therefore, 
an effective model of human cognitive performance is fundamental to the team design. %of mixed teams.
It is imperative that such a model takes into account operator's decision making mechanisms as well as  exogenous factors affecting these mechanisms, including  fatigue, situational awareness, memory, and boredom. The human operator's decisions may be erroneous, and accordingly, the quality of operator's  decision should be ascertained and incorporated in the mixed team design.
In this paper, we demonstrate the design principles for mixed human-robot teams using the context of a persistent surveillance problem.

The persistent surveillance problem involves continuous search of target regions with a team of fixed and mobile sensors.
The evidence collected by the sensor network is then processed by a human operator. 
 An efficient persistent surveillance policy has multiple objectives, including, minimizing the time between subsequent visits to a region, minimizing the detection delay at each spatial region, and maximizing visits to regions with high likelihood of an anomaly. 
%The persistent surveillance mission may be fully autonomous or may involve a human operator that processes the collected evidence and accordingly modifies the surveillance policy. 
The fundamental trade-off in persistent surveillance is between the amount of evidence collected from the visited region and the resulting delay in evidence collection from other regions. We address this trade-off by designing efficient surveillance policies that, with high probability, collect evidence from regions that are highly likely to be anomalous. We utilize cognitive models for human decision-making to ascertain the accuracy of human decisions and thus, determine the likelihood of a region being anomalous using operator's decisions. 

A second key problem is how to manage the attention of an operator in charge of analyzing the evidence collected by robotic partners. In particular, depending of the importance and the processing difficulty of the collected evidence, how much duration the operator should allocate to each piece of the collected evidence such that her overall performance is optimum. We model the attention allocation problem for human operators as a Markov decision process and  utilize certainty-equivalent receding horizon control framework to determine efficient attention allocation policies.

{\bf State of the art:} As a consequence of the growing interest in mixed teams, a significant effort has focused on modeling human cognitive performance and its integration with the automaton. Broadly speaking, there have been two approaches to the design of mixed teams in the context of surveillance. In the first approach, the human operator is allowed to take their time for decision-making on each task and the automaton is adaptively controlled to cater to the human operator's cognitive requirements. In the second approach, both the human operator and the automaton are controlled. For instance, the human operator is told the duration they should spend on each task, and their decision is utilized to adaptively control the automaton. In the first approach an  operator's performance is modeled through her reaction time on each task.  %The free reaction time is modeled as a log-normal random variable and 
The fundamental research questions under this approach include (i) optimal scheduling of the tasks to be processed by the operator~\cite{CN-BM-JWC-MLC:08, LFB-NP-MLC:10, LFB-NWMB-MLC:10, KS-CN-TT-EF:08, KS-TT-EF:08,JWC-etal:11, CEN:09}; (ii) enabling shorter operator reaction times by controlling the fraction of the total time during which  the operator is busy~\cite{KS-EF:10a, KS-EF:10b}; and (iii) efficient work-shift design to counter fatigue effects on the operator~\cite{NDP-KAM:12}.  
In the second approach an operator's performance is modeled by the probability of making the correct decision conditioned on the duration operator spends on the task.  The fundamental research questions under this approach include (i) optimal duration allocation to each task~\cite{VB-RC-CL-FB:11zc, VS-FB:12n, MM-RR:13}; (ii) controlling operator utilization to enable better performance~\cite{VS-AS-FB:11z}; and (iii)  controlling the automaton to collect the relevant information~\cite{VS-FP-FB:11za, VS-KP-FB:08p, VS-AS-FB:11z}. In this paper, we focus on the latter approach, although most of the concepts can be easily extended to the former approach.% as well.

{\bf Contributions: }The objective of this paper is to illustrate the use of systems theory to design mixed human-robot teams. We survey models of human cognition and we present a methodology to incorporate these models into the design of mixed teams. We elucidate on the design principles for mixed teams engaged in surveillance missions. 

The major contributions of this work are fourfold.
First, we present models  for several aspects of human cognition and unify them to develop a single coherent model that captures relevant features of human cognition, namely, decision-making, situational awareness, fatigue, and memory. The proposed model forms the backbone of our mixed human-robot team design.

Second, we present a certainty-equivalent receding-horizon control framework to design mixed human-robot teams. Certainty-equivalent receding-horizon control is a technique to compute the approximate solutions of dynamic optimization problems. We survey the fundamentals of certainty-equivalent receding-horizon control and apply them to design mixed team surveillance missions.

Third, we survey robotic persistent  surveillance schemes and extend them to  mixed team scenarios. In the context of mixed teams, the surveillance scheme is a component of  the closed loop system in which the robotic surveillance policy affects the human performance and human performance in turn affects the robotic surveillance policy. We present an approach to simultaneously design  robotic surveillance and operator attention allocation policies. 

Finally, we survey non-Bayesian quickest change detection algorithms. We demonstrate that human cognition models  can be used to estimate the probability a human decision being correct.  
We (i) adopt the non-Bayesian framework, and (ii) treat human decisions as binary random variables to run a change detection algorithm on human decisions and ascertain the occurrence of an anomaly within a desired accuracy. The proposed methodology can be extended to non-binary human decisions.

%\vspace{-0.3in}

\section{Human Operator Performance Modeling}

%\subsection{Elements of Human Decision Making}

\subsection{Two-alternative choice task and drift-diffusion model}
A two-alternative choice task models a situation in which an operator has to decide among one of the two alternative hypotheses. Models for two-alternative choice tasks within the continuous sensory information acquisition scenarios rely on three assumptions: (a) evidence is collected over time in favor of each alternative; (b) the evidence collection process is random; and (c) a decision is made when the collected evidence is sufficient to choose one alternative over the other. Several models for two-alternative choice tasks have been proposed; see~\cite{RB-EB-etal:06} for detailed descriptions of such models. The model considered in this paper is the drift diffusion model (DDM)~\cite{RR:78, RR-GM:08}.
We choose the DDM because: (i) it is a simple and well characterized model; (ii) it captures a good amount of behavioral and neuro-scientific data; and (iii) the optimal choice of parameters in other models for two-alternative choice tasks reduces them to the DDM~\cite{RB-EB-etal:06}. 

In the DDM  the evidence aggregation in favor of an alternative is modeled as:
\begin{equation}\label{eq:ddm}
dx(t) = \mu d t + \sigma d W(t), \quad x(0)=x_0,
\end{equation}
where $\mu\in \real$ in the drift rate, $\sigma \in \real_{>0}$ is the diffusion rate, $W(\cdot)$ is the standard Weiner
 process, and $x_0\in \real$ is the initial evidence. 
For an unbiased operator, the initial evidence $x_0=0$, while for a biased operator $x_0$ captures the odds of the prior probability that the alternative hypotheses is true; in particular, $x_0= \sigma^2\log (\pi/(1-\pi))/2\mu$, where $\pi$ is the prior probability of the first alternative being true.

\begin{figure}
\centering  \scriptsize
%\subfigure[Free Response Paradigm]{
%\includegraphics[width=0.3\textwidth]{free-response} \label{fig:free-response-ddm}
%}  
\subfigure[Interrogation Paradigm]{
\includegraphics[width=0.22\textwidth]{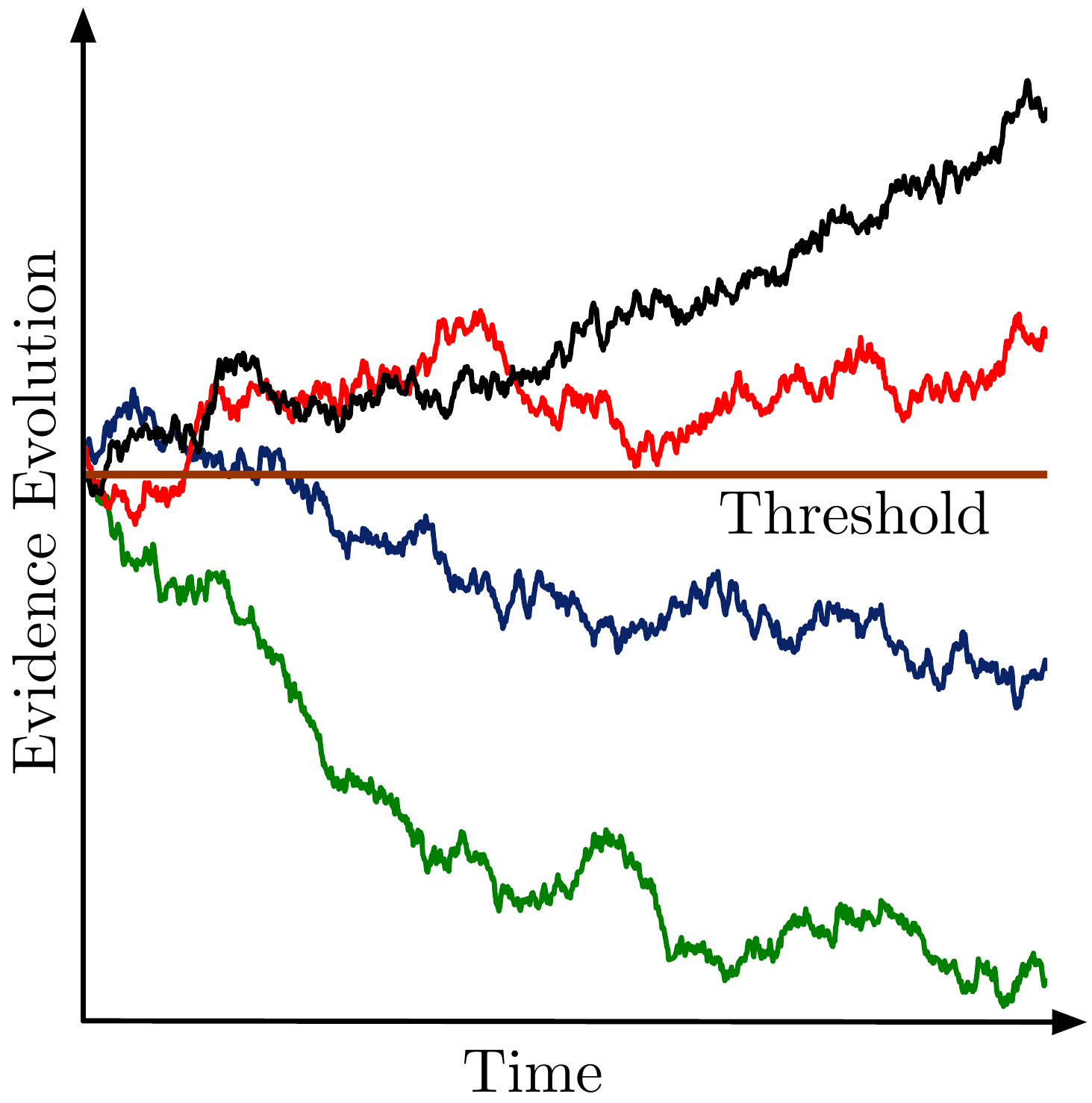} \label{fig:interrogation-ddm}
}
\subfigure[Interrogation Performance]{
\includegraphics[width=0.22\textwidth]{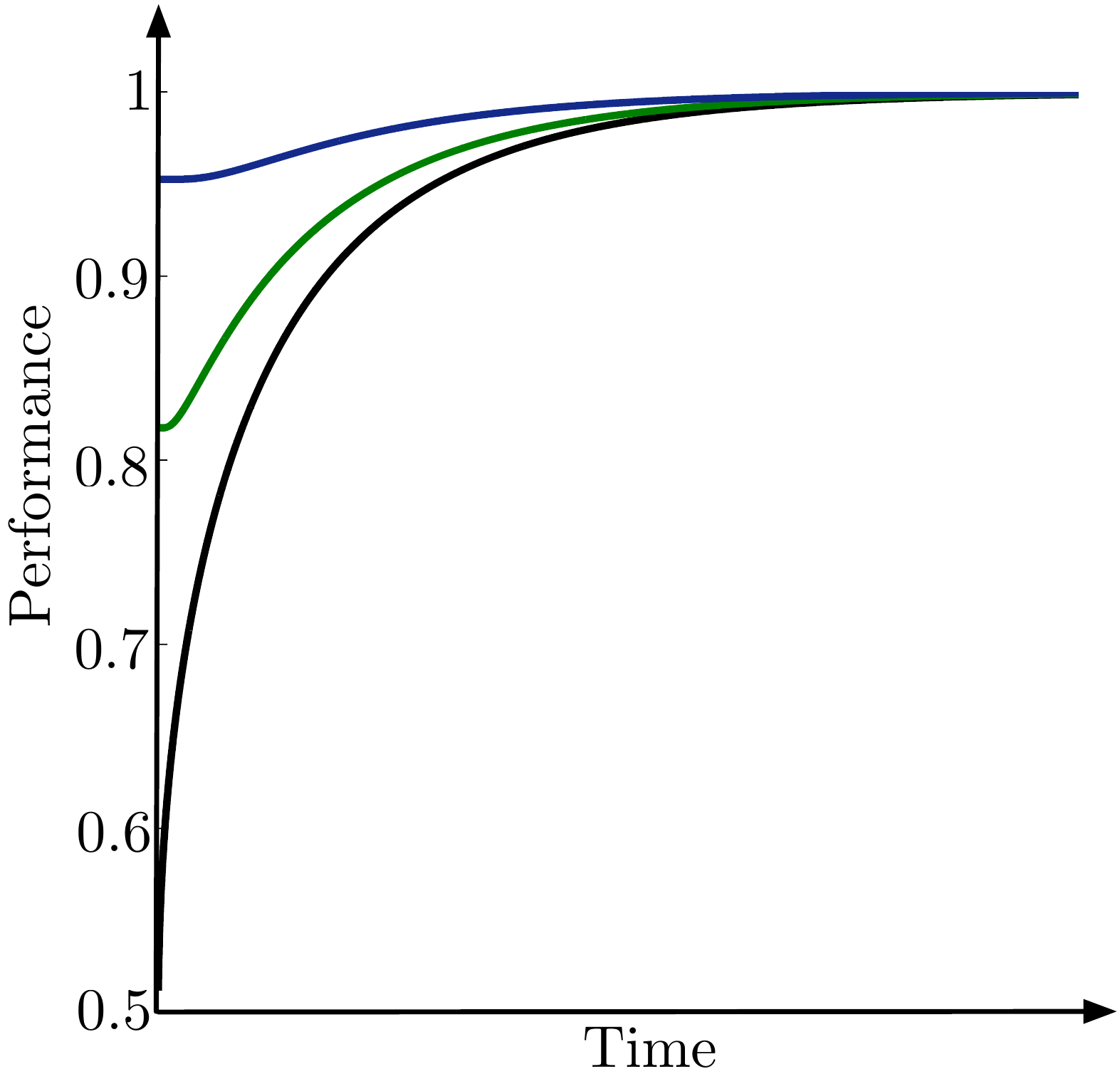}\label{fig:performance}
}
%
%\subfigure[Empirical Reaction Times]{
%\includegraphics[width=0.3\textwidth]{RT_data} \label{fig:RT_data}
%} 
%\subfigure[Reaction Times under Free Response DDM]{
%\includegraphics[width=0.3\textwidth]{RT_ddm} \label{fig:RT_ddm}
%}
\caption{\footnotesize  (a) The interrogation paradigm of decision-making. The evidence evolves according to the DDM and the decision-maker makes the correct decision if the evidence at the time of interrogation is more than a given threshold.  (b) Performance of the decision-maker under the  interrogation paradigm. The blue, green, and black curves respectively, represent performance of an operator for different initial beliefs.}
\end{figure}

For the information aggregation model~\eqref{eq:ddm}, human decision-making is studied in two paradigms, namely, \emph{free response} and \emph{interrogation}~\cite{RB-EB-etal:06}. In the free response paradigm, the operator takes her own time to decide on an alternative, while in the interrogation paradigm, the operator works under time pressure and needs to decide within a given time. 
In this section, we focus on the interrogation paradigm.
A typical evolution of the DDM under the interrogation paradigm is shown in Figure~\ref{fig:interrogation-ddm}.
The interrogation paradigm relies upon a single threshold: for a given deadline $T\in \real_{>0}$, the operator decides in favor of  the first (second) alternative if the evidence collected until time $T$, i.e., $x(T)$ is greater (smaller) than a threshold $\nu \in \real$. If the two alternatives are equally likely, then the threshold $\nu$ is chosen to be zero.
According to equation~\eqref{eq:ddm}, the evidence collected until time $T$ is a Gaussian random variable with mean $\mu T + x_0$ and variance $\sigma^2 T$. Thus, the probability to decide in favor of the first alternative is 
\[
\prob (x(T)>\nu) = 1- \prob (x(T)<\nu) = 1- \Phi\Big(\frac{\nu - \mu T-x_0}{\sigma \sqrt{T}}\Big),
\]
where $\Phi(\cdot)$ is the standard normal  cumulative distribution function.

In this paper, we use the accuracy of the decisions made by the operator as a measure of her performance.
Accordingly,  we pick the probability of making the correct decision as the performance metric. We assume the drift rates to be symmetric, i.e., the drift rates are $\pm \mu$ when the alternative $0$ and $1$ are true, respectively. Consequently, the performance function when alternative $0$ is true is $\map{f^0}{\real_{\ge 0}\times [0,1]}{[0,1)}$ defined by
\[
f^0(t, \pi)=1- \Phi\Big(\frac{\nu - \mu t-x_0}{\sigma \sqrt{t}}\Big). 
\]
Similarly, the performance function when alternative $1$ is true is $\map{f^1}{\real_{\ge 0}\times [0,1]}{[0,1)}$ defined by
\[
f^1(t, \pi)= \Phi\Big(\frac{\nu + \mu t-x_0}{\sigma \sqrt{t}}\Big). 
\]
Given the prior probability $\pi$ of the first alternative being true, the performance function $\map{f}{\real_{\ge 0} \times [0,1]}{[0,1)}$ is defined by
\begin{equation}\label{eq:performance-func}
f(t, \pi ) = \pi f^0(t, \pi) + (1-\pi) f^1(t, \pi).
\end{equation}
A typical evolution of the performance function $f$, for different values of $\pi$, is shown in Figure~\ref{fig:performance}. Note that the performance function is a sigmoid function of time.

\section{Mixed Human-Robot Team Surveillance Setup}

\begin{figure*}[ht]
\centering 
\includegraphics[width=\linewidth]{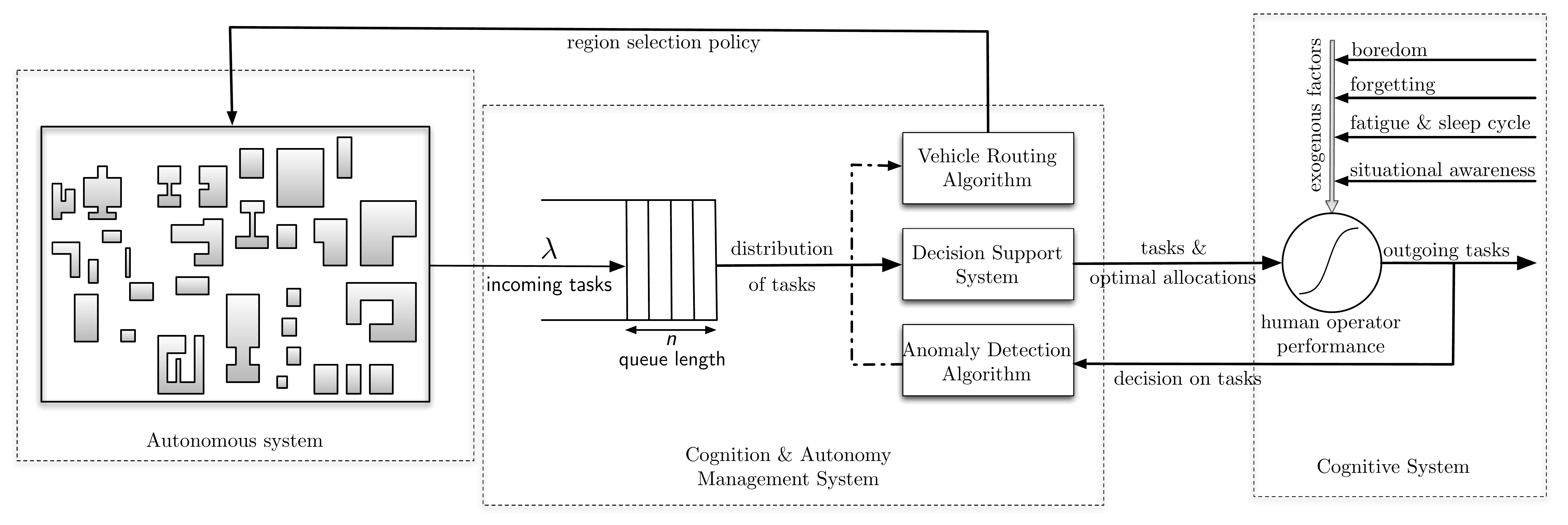}
\caption{Mixed human-robot team surveillance setup. The setup comprises three systems: (i) the autonomous system; (ii) the cognitive system; and (iii) the cognition and autonomy management system (CAMS). The autonomous system comprises a vehicle collecting information from a set of regions. The collected information is sent to the CAMS. The cognitive system comprises a human operator who looks at the evidence collected by the autonomous vehicle and decides on the presence or absence of an anomaly. The operator's performance depends on several factors including, her situational awareness, her fatigue, and her memory retention capacity. The CAMS comprises three parts (i) a vehicle routing algorithm, (ii) a decision support system, and (ii) an anomaly detection algorithm. The decision support system suggests to the operator the duration they should spend on each task. The anomaly detection algorithm treats human decision as a binary random variable and sequentially detects an anomaly within a desired accuracy.} \label{fig:problem-setup}
\end{figure*}

We study a human-in-the-loop persistence surveillance problem. The objective of the surveillance mission is to detect within a prescribed accuracy any anomaly in a set of regions. Our setup for the human-in-the-loop surveillance is shown in Figure~\ref{fig:problem-setup}. The setup comprises three components, (i) the autonomous system, (ii) the cognitive system, and (iii) the cognition and autonomy management system (CAMS). 

The autonomous system comprises  a UAV that survey a set of regions according to some routing policy. The UAV during each visit collects evidence from the region surveyed and sends it to the CAMS that eventually sends it to the cognitive component.  For simplicity, we consider only one vehicle, but our approach extends easily to multiple vehicles. 

The cognitive component comprises of a human operator who examines the evidence and decides on the absence/presence of an anomaly at the associated region. The performance of the operator depends on several variables, including, the nature of the task, her situational awareness, her fatigue, her prior knowledge about the task, and her boredom.
 
The operator sends her decision to the CAMS. The CAMS comprises of three elements: (i) the decision support system, (ii) the anomaly detection algorithm, and (iii) the vehicle routing algorithm. The decision support system based on the performance metrics of the operator suggests  the optimal duration to be allocated to a given task. The decisions made by the human operator may be erroneous. The anomaly detection algorithm is a sequential statistical algorithm that treats the operator's decision as a binary random variable and ascertains the desired accuracy of the anomaly detection. 
The anomaly detection algorithm also provides the likelihood of an anomaly at each region. The vehicle routing algorithm utilizes the likelihood of each region being anomalous to determine an efficient vehicle routing policy. 

We denote the $k$-th region by $\mc R_k, k \in \until{m}$. We model the surveillance mission as a sequence of  two-alternative choice tasks and accordingly, model the operator performance as in equation~\eqref{eq:performance-func}. 
The two alternatives in this setting are the presence of an anomaly and the absence of an anomaly, respectively.
We denote the performance of the operator at region $\mc R_k$ by 
$\map{f_k}{\real_{\ge 0} \times (0,1)}{[0,1)}$. 
We study the human-robot team surveillance under the following assumptions:
\begin{enumerate}
%\item the performance function of the operator at region $k\in \until{m}$ in absence and presence of some anomaly be  and $\map{f^1_k}{\real_{\ge 0} \times (0,1)}{[0,1]}$, respectively;
\item evidences in different regions are mutually independent;
\item the travel time between region $\mc R_i$ and $\mc R_j$ is $d_{ij}, i,j\in \until{m}$; and
\item the UAV takes time $T_k, k\in \until{m}$ to collect a sample from region $\mc R_k$.
\end{enumerate}
Queuing theory has emerged as a popular paradigm to model single-operator-multiple-vehicles systems~\cite{CN-BM-JWC-MLC:08, KS-CN-TT-EF:08, KS-TT-EF:08, KS-EF:10b, KS-EF:10a, NDP-KAM:12, VB-RC-CL-FB:11zc}. In the same vein, we model the evidence aggregated  across regions over time as  a queue of decision-making tasks. The decision-making tasks arrive according to some stochastic process and are stacked in a queue. We assume that each decision-making task arrives with a processing deadline and incorporate it as a soft constraint, namely, latency penalty. The latency penalty is the penalty incurred due to the delay in processing a given task. The operator receives a reward for a correct decision on the task. The performance function of the operator is a measure of the expected reward they obtain for allocating a given duration to a given task. The decision support system suggests to the operator duration allocations to tasks such that her overall benefit per unit task, i.e., the reward per unit task minus the latency penalty per unit task, is maximized. 

We adopt the cumulative sum (CUSUM) algorithm~\cite{HVP-OH:08} as our anomaly detection algorithm.
The details of the CUSUM algorithm are presented in the insert entitled "Quickest Change Detection."
In particular, we run $n$ parallel CUSUM algorithms (one for each region) to detect any anomaly at any of the regions. We refer to such a set of parallel CUSUM algorithms as the \emph{ensemble CUSUM algorithm}. 
\newcounter{eqncnt}
\newcounter{sideeqn}
\setcounter{eqncnt}{\value{equation}}

\setcounter{equation}{0}
\renewcommand{\theequation}{S\arabic{equation}}

\begin{figure}[t!]
\begin{shaded}
\section*{Quickest Change Detection}
Consider a sequence of observations $\{y_1,y_2,\ldots\}$ such that $\{y_1,\ldots,y_{\upsilon-1}\}$ are i.i.d. with probability density function $p^0$ and $\{y_\upsilon, y_{\upsilon +1},\ldots\}$ are i.i.d. with probability density function $p^1$ with $\upsilon$ unknown. 
The objective of the quickest change detection problem is to minimize the expected delay in the detection of the change in the underlying distribution of observations subject to a bound on the expected time between two consecutive false alarms. Specifically,  
let $\subscr{\upsilon}{det} \ge \upsilon$ be the iteration at which the change is detected. Then, the non-Bayesian quickest detection problem~\cite{HVP-OH:08, DS:85} is defined as
\begin{align}\label{eq:non-bayesian}
\begin{split}
\underset{\subscr{v}{det}}{\minimize} & \quad  \sup_{\upsilon \ge 1}  \expt_{p^\upsilon} [\subscr{\upsilon}{det} -\upsilon +1 \;|\; \subscr{\upsilon}{det}\ge \upsilon]\\
\subject & \quad \expt_{p^0}[\subscr{\upsilon}{det}] \ge \frac{1}{\gamma},
\end{split}
\end{align}
where $p^{\upsilon} \in\{p^0,p^1\}$ is the true distribution of the observation at iteration $\upsilon$, $\expt_{p^{\upsilon}}[\cdot]$ represents the expected value with respect to $p^\upsilon$,  $\gamma\in\real_{>0}$ is a small constant and is called the {\it false alarm rate.}

The solution to the problem~\eqref{eq:non-bayesian} is computed by the
cumulative sum (CUSUM) algorithm~\cite{HVP-OH:08} described as follows
\begin{algorithmic}[1]
\STATE initialize $\Lambda_0:=0$;
\STATE at time $\tau\in\naturals$, collect observation $y_\tau$;
\STATE  integrate evidence  $\Lambda_\tau:= \max\big\{0, \Lambda_{\tau-1} + \log \frac{p^1(y_\tau)}{p^0(y_\tau)}\big\}$;\\
 \emph{\% decide only if the threshold $\subscr{\Lambda}{thresh}$ is crossed}
\STATE  \textbf{if} $\Lambda_{\tau}>\subscr{\Lambda}{thresh}$, \quad \textbf{then}
  declare change is detected%
\STATE \textbf{else} ($ \Lambda_{\tau} \in {[0,\subscr{\Lambda}{thresh})}$)
  \quad continue sampling (step~2:)
\end{algorithmic}
The CUSUM algorithm has two important properties: (i) it is optimal in  that it minimizes the expected detection delay for a given false alarm rate; and (ii) it is an existing model for the decision making process for a human performing a detection task~\cite{RR-JG-MEN:03}.

If the observations $\{y_1,y_2,\ldots\}$ are dependent, then the CUSUM statistic $\Lambda_{\tau}$ is modified to 
the CUSUM-like statistic~\cite{BC-PW:00}:
\[
\Lambda_{\tau} :=  \max \Big\{0, \Lambda_{\tau-1} + \log \frac{p^1(y_\tau|y_{\tau-1}, \ldots, y_{\bar \upsilon})}{p^0(y_\tau|y_{\tau-1}, \ldots, y_{\bar \upsilon} )}\Big\},
\]
where $\bar \upsilon$ is the latest iteration at which $\Lambda_{\bar \upsilon-1}=0$.
\end{shaded}
\vspace{-0.2in}
\end{figure}

\setcounter{sideeqn}{\value{equation}}
\setcounter{equation}{\value{eqncnt}}
\renewcommand{\theequation}{\arabic{equation}}

We consider a simple routing policy that directs the vehicle to a randomly chosen region during each visit. The probability to select a region is chosen proportional to the likelihood of  a region being anomalous.   More sophisticated stochastic vehicle routing policies for surveillance are presented in the insert entitled "Stochastic Persistent Surveillance."

\setcounter{eqncnt}{\value{equation}}

\setcounter{equation}{\value{sideeqn}}
\renewcommand{\theequation}{S\arabic{equation}}

\begin{figure*}[t!]
\begin{shaded}
\section*{Stochastic Persistent Surveillance}
\medskip
\begin{minipage}{0.485\textwidth}
Persistent surveillance is a multi-vehicle mission in which a set of
autonomous vehicles endlessly visit a set of regions to collect
evidence. The collected evidence may be sent to a control center or may be
communicated among the vehicles. This evidence is utilized to detect undesired events at the set of regions under consideration. 
The vehicles visit regions according to some policy
that might be deterministic or stochastic. Deterministic policies
usually rely on (i) computing a shortest tour through the regions and (ii)
requiring the vehicles persistently move along such a tour. 
Deterministic policies are mostly periodic and predictable. 
If the event to
be detected is, for example, the existence of an intruder, then the intruder at
a region may hide only at the instants when the vehicle visits the region
and thus, the deterministic policy fails.  
%Moreover, if the persistent
%surveillance mission is characterized by a surveillance criterion that
%describes the relative number of visits that should be made to a given
%region, then a deterministic surveillance policy may not achieve the
%surveillance criterion under all conditions~\cite{KS-DMS-MWS:09}. 
% Such
%shortcomings of the deterministic policies have prompted the development of
%stochastic surveillance policies. 
Stochastic surveillance policies
typically involve an ergodic Markov chain in which each region represents a
state and transition probabilities and the stationary distribution are
designed according to an appropriate surveillance criterion.

For a single vehicle, there are two popular schemes to construct a Markov
chain with a desired stationary distribution (surveillance criterion), namely, the
Metropolis-Hastings algorithm and the fastest mixing Markov chain (FMMC)
method.  We now describe these two schemes.  Consider a set of regions
modeled by the graph $\mc G=(V, \mc E)$, where $V$ is the set of $m$ nodes
(each node corresponds to a region) and $\mc E$ is the set of edges
representing the connectivity of the regions.  Let the surveillance
criterion be $\q = (q_1,\dots,q_m) \in \Delta_{m}$.  The
Metropolis-Hastings algorithm~\cite{LW:04} picks the transition matrix $A$,
i.e., the matrix of transition probabilities from each state to every other
state of the Markov chain, according to:
\[
A_{ij}=\begin{cases}
  0, & \text{if } (i,j) \notin \mc E,\\
  \min \big\{\frac{1}{d_i}, \frac{q_j}{q_i d_j}\big\}, & \text{if } (i,j) \in \mc E \text{ and } i\ne j,\\
  1-\sum_{k=1, k\ne i}^n A_{ik}, & \text{if } (i,j) \in \mc E \text{
    and } i= j,
\end{cases}
\] 
where $d_i$ is the number of regions that can be visited from region $\mc R_i$.
%\smallskip
\end{minipage}
\hfill
\begin{minipage}{0.485\textwidth} 
The transition matrix $A\in
\real^{m\times m}$ of the FMMC  with a desired stationary distribution $\q \in
\Delta_{m}$ is determined by solving the following semi-definite program~\cite{SB-PD-LX:04}:
\begin{equation*}
\begin{split}
\minimize & \quad \| Q^{1/2} A Q^{1/2} -\subscr{\q}{root} \subscr{\q \tran}{root} \|_2\\
\subject & \quad A \boldsymbol{1}_m =\boldsymbol{1}_m\\
& \quad QA=A{\tran} Q \\
&\quad A_{ij}\ge 0, \text{ for each } (i,j) \in \mc E \\
&\quad A_{ij}=0, \text{ for each } (i,j)\notin \mc E,
\end{split}
\end{equation*}
where $Q$ is a diagonal matrix with diagonal $\q$,
$\subscr{\q}{root}=(\sqrt{q_1},\ldots, \sqrt{q_n})$, and $\boldsymbol{1}_m$
is the vector of all ones.  In order to achieve the surveillance criterion
at an accelerated rate, a time-varying Markov chain can also be constructed in
the spirit of~\cite{JG-JB:05}.

For multiple vehicles, a naive stochastic policy that achieves the
surveillance criterion is to let each vehicle follow the single vehicle
policy. One potential drawback of such a naive policy is that two or more
vehicles may survey the same region at the same time. This drawback can be
removed by constructing a Markov chain on a lifted space from which the
undesired states are removed.  Such a policy is developed
in~\cite{KS-DMS-MWS:09} and a message passing based auction algorithm is
used to design a decentralized policy that achieves the desired
surveillance criterion.

Another important feature of the stochastic surveillance is the choice of the surveillance criterion. In this paper, we choose the surveillance criterion proportional to the likelihood of the anomaly at each region. 
In general, the surveillance criterion depends on the mission objective. For example,
if the mission objective is the detection of anomalous regions, then the surveillance criterion may be chosen to minimize the average detection delay~\cite{VS-FP-FB:11za}. The minimization of the average detection delay inherently considers the difficulty of detection at each region, the travel times between the regions, and the likelihood of each region being anomalous.  

\medskip
\end{minipage}
\end{shaded}
\end{figure*}

\setcounter{sideeqn}{\value{equation}}
\setcounter{equation}{\value{eqncnt}}
\renewcommand{\theequation}{\arabic{equation}}

%\subsubsection{A unified model for human operator}

\section{Design of the Decision Support System}

%\subsection{Design under Simplified Operator Model}
We first consider the design on the decision support system. In this section, we ignore the exogenous factors and their effects on the human performance. 
We design the  decision support system  under the following assumptions:
\begin{enumerate}
\item performance functions of the operator on a task from region $\mc R_k$ in absence and presence of an anomaly are $\map{f^0_k}{\real_{\ge 0} \times (0,1)}{\real_{\ge 0}}$ and $\map{f^1_k}{\real_{\ge 0} \times (0,1)}{\real_{\ge 0}}$, respectively;
\item  each task from region $\mc R_k$ comes with a processing deadline $\supscr{T}{ddln}_k \in \real_{\ge 0}$;
\item based on the importance of the region, a weight $w_k \in \real_{>0}$ is assigned to each task collected from region $\mc R_k$;  
\item tasks arriving to the queue while the $\ell$-th task is served are sampled from a probability distribution that assigns a probability $q_k^\ell \in (0,1)$ to region $\mc R_k$.  
\end{enumerate}
Similar to~\eqref{eq:performance-func}, the average performance function $\map{f_k}{\real_{\ge 0} \times (0,1)}{\real_{\ge 0}}$ at region $\mc R_k$ is defined by
\[
f_k(t, \pi)= (1-\pi)f_k^0(t, \pi) + \pi f_k^1 (t, \pi). 
\]
Under the aforementioned assumptions, each task from region $\mc R_k$ is characterized by the triplet $(f_k,  w_k, \supscr{T}{ddln}_k)$. 

Decision-making tasks collected by the vehicle are stacked in a queue and are processed by the operator in a first-come-first-serve processing discipline. The task arrival process for this queue depends on the vehicle routing policy.
Let the $\ell$-th task in the queue be from region  $\mc R_{k_\ell}$ and let $n_\ell$ be the queue length
when the processing of the $\ell$-th task is initiated. Let the belief of the operator about region $\mc R_{k}$ being anomalous before processing the $\ell$-th task be $\pi_{k}^{\ell-1}$. Without loss of generality, we assume that initially the operator is unbiased about each region being anomalous, i.e., $\pi_{k}^{0}=0.5$, for each $k\in \until{m}$. Given a duration allocation $t_\ell \in \real_{>0}$ to the  $\ell$-th task in the queue, the operator's belief after processing the $\ell$-th task  can be estimated using the Bayes rule as follows:
\[
\bar \pi_{j}^{\ell} =
\begin{cases}
 \frac{ \pi_{j}^{\ell-1} \prob({\texttt{dec}}_\ell| H^1_{k}, t_\ell)}{(1- \pi_{j}^{\ell-1}) \prob({\texttt{dec}}_\ell| H^0_{k_\ell}, t_\ell)+ \pi_{j}^{\ell-1} \prob({\texttt{dec}}_\ell| H^1_{k_\ell}, t_\ell)}, & \text{if } j=k_{\ell}, \\
 \pi_{j}^{\ell-1}, & \text{otherwise,}
 \end{cases}
\]
where $H^0_k$ and $H^1_k$ denote the hypothesis that region $\mc R_k$ is non-anomalous and anomalous, respectively, and ${\texttt{dec}}_\ell \in \{0,1\}$ is the operator's decision,
%\[
%\bar \pi_{k} = \frac{ \exp \big \{\log\big( \frac{\supscr{\pi_{k}}{last}}{1- \supscr{\pi_{k}}{last}}\big) \texttt{rem}(\supscr{T}{last}_{k} )\big \}}
%{1+ \exp \big \{\log\big( \frac{\supscr{\pi_{k}}{last}}{1- \supscr{\pi_{k}}{last}}\big) \texttt{rem}(\supscr{T}{last}_{k} )\big \}}, 
%\]
and $\prob({\texttt{dec}}_\ell|\cdot, t_\ell)$ is determined from the performance function of the operator, i.e., 
\[
\prob({\texttt{dec}}_\ell=1| H^1_{k_\ell}, t_\ell) = f^1_{k_\ell} (t_\ell, \pi_{k_\ell}^{\ell-1}).
\]
%\[
%\bar \pi_{j}^{\ell} = 
%\begin{cases}
%\frac{\pi_{j}^{\ell-1} f^1_{j}(t_\ell, \pi_{j}^{\ell-1}) }
%{(1-\pi_{j}^{\ell-1}) f^0_{j}(t_\ell, \pi_{j}^{\ell-1})+ \pi_{j}^{\ell-1} f^1_{j}(t_\ell, \pi_{j}^{\ell-1})}, & \text{if } j=k_\ell, \\
%\pi_{j}^{\ell-1}, & \text{otherwise}.
%\end{cases}
%\]
In a surveillance mission, the event that a region becomes anomalous corresponds to a change in the 
characteristic environment of the region and this change may happen at an arbitrary time. It is evident from the neuroscience literature~\cite{RR-JG-MEN:03} that in a sequential change detection task, if the belief of the operator about a region being anomalous is below a threshold, then the operator resets her belief to the threshold value. Without loss of generality, we choose this threshold as $0.5$. Consequently, if the belief of the operator about a region being anomalous is below $0.5$, then they reset it to $0.5$, i.e., the belief of the operator at region $\mc R_k$ after processing the $\ell$-th task is
\[
\pi^\ell_k =\max\{0.5,\bar \pi_{j}^{\ell} \}.
\]
We incorporate the deadline on a task as a soft constraint (latency penalty) and assume that the task from region $\mc R_k$ loses value while waiting in the queue. Let the task from region $\mc R_k$ loses value at a rate $c_k$ per unit delay in its processing. The latency penalty rate $c_k$ is a function of task parameters $(f_k,  w_k, \supscr{T}{ddln}_k)$. For simplicity of notation, we drop the arguments of $c_k$. The performance function of the operator depends on her belief about regions being anomalous which varies over the evolution of the surveillance mission. Consequently, the performance function of the operator is time-varying and hence, the latency penalty rate $c_k$ is also time-varying. We denote the latency penalty rate at region $\mc R_k$ when the $\ell$-th task is processed by $c_k^{\ell}$.

%We incorporate the processing deadline on a task from region $\mc R_k$ 

The decision support system suggests to the operator the duration allocation to each task. To this end, the decision support system maximizes the infinite-horizon average reward for the human operator. 
The reward  $\map{r}{\until{m}^{n_\ell}\times (0,1)^m\times \real_{\ge 0}}{\real}$ obtained by allocating duration $t$ to the $\ell$-th task is defined by

\begin{multline*}
r (\bs k_\ell, \bs \pi^{\ell-1}, t) = w_{k_\ell}f_{k_\ell}(t,  \pi_{k_\ell}^{\ell-1}) \\ 
- \frac{1}{2}\Big( \sum_{i=\ell}^{\ell+n_\ell-1} c_{k_i}^{\ell} + \sum_{j=\ell}^{\ell+n'_\ell-1} c_{k_j}^{\ell} \Big) t,
%- \frac{1}{2} \bar{c} \lambda t^2,
\end{multline*}
where  $\boldsymbol k_\ell \in \until{m}^{n'_\ell}$ is the vector of region indices of the tasks in the queue,  $\bs \pi^{\ell-1}$ is the vector of operator's belief about each region being anomalous, and $n'_\ell$ is the queue length just before the end of processing the $\ell$-th task. 
Note that the queue length while a task 
is processed may not be constant, therefore, the latency penalty is computed as the average of the latency penalty for the tasks present
at the start of processing the task and  the latency penalty for the tasks present at the end of processing the task.  
Let the decision support system suggest to the operator to allocate duration $t_\ell$ to the $\ell$-th task in the queue. Accordingly, the objective of the decision support system is to maximize the infinite-horizon average reward defined by
\[
V_{\text{avg}} = \liminf_{N\to+\infty} \frac{1}{N} \sum_{\ell=1}^N  \expt[r (\boldsymbol k_\ell, \bs \pi^{\ell-1}, t_\ell)].
\]

\setcounter{eqncnt}{\value{equation}}

\setcounter{equation}{\value{sideeqn}}
\renewcommand{\theequation}{S\arabic{equation}}

\begin{figure}[t!]
\begin{shaded}

\section*{Certainty-Equivalent Receding-Horizon Control} 
%If the optimal value function for the MDP exists, then it needs to be efficiently computed. The exact computation of the optimal value function is in general intractable and approximation techniques are 
%employed~\cite{DB:05}. 
Consider a time-varying dynamical system of the form
\[
x_{\ell+1} = \texttt{evol}_\ell(x_\ell, v_\ell, d_\ell),
\]
where $\ell\in \naturals$ denotes the discrete time, $x_\ell \in  \mc X$ is the state variable, $v_\ell \in  \mc V$ is the control signal,  $d_\ell \in  \mc D$ is the disturbance and $\map{\texttt{evol}_\ell}{\mc X \times \mc V \times \mc D}{\mc X}$ denotes the evolution map of the system. 
Given a stochastic model for $\seqdef{d_\ell}{\ell\in\naturals}$,
the control objective is to design control sequence $\bs v:=\seqdef{v_\ell}{\ell \in \naturals}$ that minimizes:
\begin{equation}\label{eq:infinite-horizon}
\subscr{J}{avg} (x_1, \bs v) = \limsup_{N\to +\infty} \frac{1}{N} \sum_{\ell=1}^N \expt [g(x_\ell,  v_\ell, {d}_{\ell})],
\end{equation}
where $x_1$ is the initial state, $\bs v$ is the control sequence, 
$\map{g}{\mc X \times \mc V \times \mc D}{\real}$ is the stage cost,
and the expectation is taken over the disturbance $d_\ell$.

The minimization of $\subscr{J}{avg}$ is a hard optimization problem and is usually approximately solved using receding-horizon control.
The receding-horizon control approximates the infinite-horizon average value function with a finite horizon average value function and solves such a finite horizon problem at each iteration to determine the control at that iteration.  We focus on a particular type of receding-horizon control, namely, 
the certainty-equivalent receding-horizon control~\cite{DB:05,HSC-SIM:03, JM-YW-SB:11}. 
 According to the
certainty-equivalent approximation, the future uncertainties in the finite horizon evolution are replaced by their expected values.
Hence, the certainty-equivalent receding-horizon control solves the following optimization problem to determine the control to apply at iteration $\ell$ 
 \begin{equation}\label{eq:receding}
\underset{\hat v_1, \ldots, \hat v_N}{\maximize} \quad \frac{1}{N} \sum_{j=0}^{N-1} g(\bar x_{\ell+j},  \hat v_j,  \bar d_{\ell+j}), 
 \end{equation}
where $\{\hat v_1, \ldots, \hat v_N\}$ is the control sequence, and $\bar x_{\ell+j}$ is the certainty-equivalent evolution of the system, i.e., the evolution of the system obtained by replacing the uncertainty in the evolution at each stage by its expected value, $\bar x_\ell =x_\ell$, and $\bar d_{\ell+j}$ is the expected value of the uncertainty at stage $\ell+j$. 
% the approximate average value function $\map{\subscr{\bar J}{avg}}{\mc X}{\real}$ is defined by
%\[
%\subscr{\bar J}{avg} (x) = \max_{u_1, \ldots, u_N} \frac{1}{N} \sum_{\ell=1}^N g(\bar x_\ell,  u_\ell, \bar{w}_{\ell}).
%\]
%where  
The certainty-equivalent receding-horizon control at iteration $\ell$ solves optimization problem~\eqref{eq:receding} and picks $v_\ell = \hat v_1$. 
%computes the approximate average value function $\subscr{\bar J}{avg} (x) $, implements the first control action $u_1$, lets the system evolve to a new state $x'$, and repeats the procedure at the new state.
There are two salient features of the certainty-equivalent receding-horizon control: (i) it approximates the value function at a given state by a deterministic dynamic program; and (ii) it utilizes a open loop strategy, i.e., the action variable $\hat v_j$ does not depend on $\bar x_{\ell+j}$. Thus, if the deterministic dynamic program~\eqref{eq:receding} can be solved efficiently, then the certainty-equivalent receding-horizon control offers a computationally tractable sub-optimal solution to problem~\eqref{eq:infinite-horizon}.
\end{shaded}
\end{figure}
%\vspace{-0.3in}

\setcounter{sideeqn}{\value{equation}}
\setcounter{equation}{\value{eqncnt}}
\renewcommand{\theequation}{\arabic{equation}}

We adopt the certainty-equivalent receding-horizon control~\cite{DB:05, HSC-SIM:03, JM-YW-SB:11} as the algorithm for the decision support system to approximately compute the policy that maximizes the infinite-horizon average reward, and hence, obtain efficient allocations for the operator. Fundamentals of certainty-equivalent receding-horizon control are presented in the insert entitled "Certainty-Equivalent Receding-Horizon Control."

Under the certainty-equivalent approximation, the future uncertainties of the system are approximated by their expected values~\cite{DB:05}. The expected value of the belief over the prediction horizon is equal to the current belief, i.e., $\expt[\bar \pi_k^{\ell}|\pi_k^{\ell-1}] =\pi_k^{\ell-1}$. 
Accordingly, the  travel time and evidence aggregation time under the certainty-equivalent approximation while the $\ell$-th task  is processed is ${\bs q^\ell} \tran D {\bs q^\ell} + {\bs q^\ell}\tran \bs T$. 
Furthermore, the evolution of the queue length under the certainty-equivalent approximation while the $\ell$-th task is processed is
\[
\bar{n}_{\ell+j+1} = \max \{ 1, \bar{n}_{\ell+j} -1 + {\lambda_\ell \bar t_{j}}\},
\]
where $\bar{n}_{\ell+j}$ represents predicted queue length  at the start of processing the $(\ell+j)$-th task,
$\bar t_{j}$ is the duration allocation to the $(\ell+j)$-th task,  $\bar n_\ell = n_\ell$, and $\lambda_\ell =1 /{({\bs q^\ell} \tran D {\bs q^\ell} + {\bs q^\ell}\tran \bs T)}$. 

Moreover, under the certainty-equivalent approximation, the parameters of the tasks that have not yet arrived are replaced by their expected values. While the $\ell$-th task is processed, the expected value of the performance function  $\map{\bar{f}_\ell}{\real_{\ge 0}\times [0,1]}{{[0,1]}}$, the expected importance $\bar w^\ell$ and the expected latency penalty rate $\bar{c}^\ell$ are defined by
%Accordingly, under the certainty-equivalent approximation while the $\ell$-th task is processed, the performance function, the importance, and the latency penalty rate for the tasks that have not yet arrived are, respectively, replaced by the expected performance function  $\map{\bar{f}_\ell}{\real_{\ge 0}\times [0,1]}{{[0,1]}}$, the expected importance $\bar w^\ell$, and the expected latency penalty rate $\bar{c}^\ell$ defined by
%
\begin{align*}
\fav_\ell(t, \pi_{k}^{\ell-1})& = \frac{\sum_{k=1}^m q_k^\ell w_k f_k(t, \pi_{k}^{\ell-1})}{\sum_{k=1}^m q_k^\ell w_k}, \\
  \bar{w}^\ell & =   \sum_{k=1}^m q_k^\ell w_k, \quad 
\text{and} \quad \bar{c}^\ell = \sum_{k=1}^m q_k^\ell c_k^\ell,
\end{align*}
respectively.

In the receding-horizon framework~\cite{JM-YW-SB:11}, a finite-horizon optimization problem at each iteration is solved. 
Consider a horizon length $N$, the current queue length $n_\ell$,  realizations of the sigmoid functions associated with the tasks in the queue $f_{k_{\ell+j-1}}, j\in \until{n_\ell}$,
the associated latency penalties $c_{k_{\ell+j-1}}^\ell, j\in \until{n_\ell}$ and weights $w_{k_{\ell+j-1}}, j\in \until{n_\ell}$. If $n_\ell < N$, then we define
the reward $\map{r_j}{\real_{\ge 0}}{\real}$, $j\in\until{N}$, associated with the $(\ell+j-1)$-th task by 
\begin{equation*}%\label{eq:reward-str}
  r_j(t_j) =\begin{cases}
  r^{\text{rlzd}}_j(t_j), \quad & \text{if }1\le  j \le n_\ell,\\
 r^{\text{exp}}_j(t_j), & \text{if } n_\ell+1 \le j \le N,
  \end{cases}
 \end{equation*}
where $ r^{\text{rlzd}}_j$ is the reward associated with a task whose parameters have been realized and is defined by
\begin{multline*}
  r^{\text{rlzd}}_j(\bar t_j) =w_j f_{k_{\ell+j-1}}(\bar t_j) - \frac{1}{2}\bar{c}^\ell \lambda_\ell \bar t_j^2 \\
  - \Big(\sum_{i=j}^{n_\ell} c_{k_{\ell+j-1}}^\ell +(\bar n_{\ell+j-1} - n_\ell -j+1) \bar{c} ^\ell \Big) \bar t_j , 
\end{multline*} 
and $r^{\text{exp}}_j$ is the reward associated with a future task whose parameters are approximated by their expected value and is defined by
\[ 
r^{\text{exp}}_j(\bar t_j) = \bar{w}^\ell \fav_\ell (\bar t_j) -\bar{c}^\ell \bar n_{\ell+j-1} \bar t_j - \frac{1}{2}\bar{c}^\ell\lambda_\ell \bar t_j^2.
\]
Similarly, if $n_\ell\ge N$, then the parameters associated with each task in the prediction horizon have been realized and consequently, the reward is defined by $r_j=  r^{\text{rlzd}}_j$, for each  $j\in\until{N}$. 

It should be noted that under the the certainty-equivalent approximation, the performance functions and latency penalties are time-invariant over the prediction horizon. In particular, the time-varying component of the performance function is the belief of the operator and the expected value of the belief over the prediction horizon is equal to the current belief, i.e., $\expt[\bar \pi_k^{\ell}|\pi_k^{\ell-1}] =\pi_k^{\ell-1}$.   Similarly, the latency penalty associated with each region is also a constant over the prediction horizon.

The certainty-equivalent receding-horizon policy solves the following finite-horizon optimization problem to determine the allocation to the $\ell$-th task:
\begin{equation}\label{eq:maximize-receding-horizon-real-time}
\begin{split}
\underset{\bar{\t} \succeq0}{\maximize}& \quad \frac{1}{N}\sum_{j=1}^{N} r_j (\bar t_j) \\ 
\subject&\quad \bar n_{\ell+ j+1}=  \max \{ 1, \bar{n}_{\ell+ j} -1 +  \lambda_\ell \bar t_{j}\},
\end{split}
\end{equation}
where $\bar{n}_\ell=n_\ell$ and $\bar{\t}=\{\bar t_1,\dots, \bar t_N\}$ is the duration allocation vector. In particular, if $\{\bar t_1^*,\dots, \bar t_N^*\}$ is the solution to~\eqref{eq:maximize-receding-horizon-real-time}, then the  certainty-equivalent receding-horizon policy allocates a duration $t_\ell = \bar t_1^*$ to the $\ell$-th task. 

The optimization problem~\eqref{eq:maximize-receding-horizon-real-time} is a finite-horizon dynamic program with univariate state and control variables. 
Moreover, the reward function and the state evolution function in optimization problem~\eqref{eq:maximize-receding-horizon-real-time} are Lipschitz continuous. The control variable on a task from region $\mc R_k$ is upper bounded by the deadline of the task. The certainty-equivalent queue length can be bounded above by a large constant (see~\cite{VB-RC-CL-FB:11zc} for a formal argument). 
These properties imply that the control and the state space in problem~\eqref{eq:maximize-receding-horizon-real-time} can be efficiently discretized to determine a solution that yields a value within $\epsilon$ of the optimal value~\cite{DPB:75}. Furthermore, such a discretized dynamic program can be solved using the backward recursion algorithm~\cite{DPB:01a} in $O(N/\epsilon^2)$ time. More details on discretizing a dynamic program are presented in the insert entitled "Discretization of the Action and the State Space."
  
The above methodology for designing the decision support system assumes that latency penalty is known as a function of task parameters and the region selection policy. If the deadline on a task is more than the inflection point of the associated sigmoid performance function, then one choice of such a latency penalty  is
 $c_k^\ell=w_k f'_k(\supscr{T}{ddln}_k, \pi_k^{\ell-1})$, where $f'$ represents derivative of $f$ with respect to $t$.
 % and $\pi_k$ is the belief of the operator on the presence of anomaly at region $\mc R_k$. 
The above procedure for the design of the support system yields a non-zero allocation to each task only if the latency penalty is sufficiently small. At large latency penalties, the loss in the value of the task is more than the reward obtained by processing it, and consequently, the optimal policy allocates zero duration to each task. The latency penalty should be small enough such that if only one task is present, it should be allocated a non-zero duration. This can be achieved by choosing a large enough deadline on the task. In particular, the latency penalty at region $\mc R_k$ should be small enough such that
 \[
 \text{argmax} \setdef{w_k f_k(t, \pi_k^{\ell-1}) - c_k^\ell t -\bar{c}^\ell \lambda_\ell t^2/2}{t \in \real_{\ge 0}} >0.
 \]
 
Unfortunately, for large values of $\pi_k^{\ell-1}$, the reward rate is very low
and consequently, the admissible penalty rates are very low. 
Such low penalty rates require the deadlines to be very large, which is undesirable. In the following, we will assume that the deadline is large enough for a moderately large critical value of the belief, e.g., $\pi_k^{\ell-1}\le 0.8$, and once the operator's belief is above this critical value, then a duration equal to the deadline on the task is allocated to the task. Such an allocation will bring the operator's belief below the critical value by either (i) getting an evidence that diminishes the belief, or (ii) crossing a threshold and resetting the operator's belief (see the anomaly detection algorithm for more details).

 The formulation~\eqref{eq:maximize-receding-horizon-real-time} of the certainty-equivalent finite-horizon problem ensures that the decision-making queue is stable provided the latency penalty rate for each task $c_k$ is positive. If there is no deadline on any task, then the latency penalty on each task is zero and the formulation~\eqref{eq:maximize-receding-horizon-real-time} does not ensure the stability of the decision-making queue. However, in the case of no deadlines and for a policy that allocates a duration $\supscr{t}{reg}_k\in \real_{\ge 0}$ to a task from region $\mc R_k$,  
the  average value function under the certainty-equivalent approximation  while the $\ell$-th task is processed is $\sum_{k=1}^m q_k^{\ell} w_k f_k(\supscr{t}{reg}_k, \pi^{\ell-1}_k)$. 
Thus, the allocation to the $\ell$-th task under the certainty-equivalent receding-horizon policy is determined
by the solution of the following optimization problem:
\begin{align}\label{eq:knapsack-sigmoid}
\begin{split}
\maximize & \quad \sum_{k=1}^m q_k^{\ell} w_k f_k(\supscr{t}{reg}_k, \pi^{\ell-1}_k) \\
\subject & \quad  \sum_{k=1}^m q_k^{\ell} \supscr{t}{reg}_k \le \frac{1}{\lambda_\ell}\\
& \quad \supscr{t}{reg}_k \ge 0, \text{ for each } k\in\until{m},
\end{split}
\end{align}
where the first constraint ensures the stability of the queue. 
In particular, a duration ${\supscr{t}{reg}_{k_\ell}}^*$ is allocated to the $\ell$-th task, where $({\supscr{t}{reg}_{1}}^*,\ldots, {\supscr{t}{reg}_{m}}^*)$ is the solution to~\eqref{eq:knapsack-sigmoid}.
The optimization problem~\eqref{eq:knapsack-sigmoid} is a knapsack problem with sigmoid utilities. This problem is NP-hard and a procedure to compute a $2$-factor solution is presented in the insert entitled "Knapsack Problem with Sigmoid Utilities." 

  \setcounter{eqncnt}{\value{equation}}

\setcounter{equation}{\value{sideeqn}}
\renewcommand{\theequation}{S\arabic{equation}}

\begin{figure}[t!]
\begin{shaded}
\section*{Discretization of the Action and the State Space}
Consider a time-varying deterministic dynamical system of the form
\[
x_{\ell+1} = \texttt{evol}_\ell(x_\ell, v_\ell),
\]
where $\ell\in \naturals$ denotes the discrete time, $x_\ell \in\mc X$ is
the state variable, $v_\ell \in \mc V$, and $\map{\texttt{evol}_\ell}{\mc X
  \times \mc V \times \mc D}{\mc X}$ denotes the evolution map of the
system. Let $\mc X$ and $\mc V$ be continuous and compact spaces. The
objective is to solve the following finite-horizon optimization problem:
 \begin{equation}\label{eq:finite}
\underset{v_1, \ldots, v_N}{\maximize} \quad \frac{1}{N} \sum_{j=0}^{N-1} g(x_{\ell+j},  v_j), 
 \end{equation}
 where $\map{g}{\mc X \times \mc V}{\real}$ is the stage cost.
The optimization problem~\eqref{eq:finite} is a finite-horizon dynamic program with a compact state and action space.  A popular technique to approximately solve such dynamic programs involves discretization of action and state space followed by the backward induction algorithm. We focus on the efficiency of such discretization. 
%Let the certainty-equivalent evolution of the state be represented by 
%$\bar x_{\ell+1} = \texttt{evol}(\bar x_\ell, u_\ell, \bar{w}_\ell)$, where $\map{\texttt{evol}}
%{\bar{\mc X}_\ell \times \mc U_\ell \times \mc W_\ell}{\bar {\mc X}_{\ell+1}}$ represents the certainty-equivalent evolution function, $\bar{\mc X}_\ell$ represents the compact space to which the certainty-equivalent state at stage $\ell$ belongs, $\mc U_\ell$ is the compact action space at stage $\ell$ and $\mc W_\ell$ is the disturbance space at stage $\ell$.  
Let the action space and the state space be discretized (see~\cite{DPB:75} for details of discretization procedure) such that $\Delta x$ and $\Delta u$ are the maximum grid diameter in the discretized state and action space, respectively.  Let 
$\map{\subscr{\hat{J}}{avg}^*}{\mathcal{X}}{\real}$ 
be the approximate average value function obtained via discretized state and action space and let $\map{\subscr{\bar{J}}{avg}^*}{\mathcal{X}}{\real}$ be the optimal value of the optimization problem~\eqref{eq:finite}. 
%Consider a $N$-horizon dynamic program with continuous state and action variables that belong to compact sets. Let $\mathcal{X}_k$ is the subset of Euclidean space in which the state variable lies at stage $k$, and $\map{J^*_k}{\mathcal{X}_k}{\real}$ be the optimal value function at stage $k$. 
%\begin{proposition}[\bit{Proposition~2, }]\label{prop:discretization}
If the action and state variables are continuous and belong to  compact sets, and 
the reward function and the state evolution map are Lipschitz continuous, then 
\[
|\subscr{\bar J}{avg}^*(x) - \subscr{\hat{J}}{avg}^*(x)| \le \beta(\Delta x+\Delta u),
\] 
for each $x\in \mathcal{X}$, where $\beta$ is a constant independent of the discretization grid~\cite{DPB:75}.

%The control objective is to design control sequence $\bs v:=\seqdef{v_\ell}{\ell \in \until{N}}$ such that the following cost function is minimized
%\begin{equation}\label{eq:infinite-horizon}
%\subscr{J}{avg} (x_1, \bs v) = \limsup_{N\to +\infty} \frac{1}{N} \sum_{\ell=1}^N \expt [g(x_\ell,  v_\ell, {d}_{\ell})],
%\end{equation}
%where $x_1$ is the initial state, $\bs v$ is the control sequence, and the expectation is taken over the disturbance $d_\ell$.

%Note that the certainty-equivalent evolution $\bar x_\ell$ may not belong 
%to $\mc X$, for instance, $\mc X$ may be the set of natural numbers and $\bar x_\ell$ may be a positive real number. In general, certainty-equivalent
% state $\bar x_\ell$ will belong to a compact set and therefore, the finite-horizon  deterministic dynamic program associated with  the computation of the approximate average value function involves compact action and state space. 

%\end{proposition}
\end{shaded}
\end{figure}

\setcounter{sideeqn}{\value{equation}}
\setcounter{equation}{\value{eqncnt}}
\renewcommand{\theequation}{\arabic{equation}}

 \setcounter{eqncnt}{\value{equation}}

\setcounter{equation}{\value{sideeqn}}
\renewcommand{\theequation}{S\arabic{equation}}

\begin{figure}[t!]
\begin{shaded}

\section*{Knapsack Problem with Sigmoid Utilities}

Consider a set of items $k\in\until{N}$ such that the utility of item $k$
is a sigmoid function $\map{f_k}{\real_{\ge 0}}{\real_{\ge 0}}$ of the
resource allocated to it. Let $T\in\real_{>0}$ be the total available
resource. The knapsack problem with sigmoid utilities is to determine the
optimal resource to be allocated to each item such that the total utility
is maximized. In particular, the knapsack problem with sigmoid utilities is
\begin{equation*} %\label{eq:hybrid-knapsack}
  \begin{split}
\maximize &\quad \sum_{k=1}^{N} w_{k}f_{k}(t_{k})\\
    \subject  &\quad  \sum_{k=1}^{N} t_{k} \le T\\
              &\quad t_k \ge 0, \text{ for each } k\in\until{N},
  \end{split}
\end{equation*}
where $w_k\in\real_{\ge 0}, k\in\until{N}$ is the weight associated with item $k$. 

The knapsack problem with sigmoid utility is an NP-hard problem and a $2$-factor solution to it is proposed in~\cite{VS-FB:12n}.
Before we describe  the solution, we introduce the following notation.  
For a given sigmoid function $f_k$ with inflection point $\tinf_k$, let  $\map{\fdi_k}{\real_{> 0}}{\real_{\ge 0}}$ denote the pseudo-inverse of the derivative of $f_k$ defined by
\begin{equation*}
\fdi_k(y) =\begin{cases}
\max\setdef{t\in\real_{\ge 0}}{f'_k(t)=y}, & \text{if } y\in {(0,f'_k(\tinf_k)]}, \\
0, & \text{otherwise.}
\end{cases}
\end{equation*}
Define $\alpha_{\max}= \max\setdef{w_k f'_k(\tinf_k)}{k\in\until{N}}$. We also define $\map{\subscr{F}{LP}}{{(0,\alpha_{\max}]}}{\real_{\ge 0}}$ as the optimal value of the objective function  in the following  $\alpha$-parametrized fractional knapsack problem:
\begin{align*}%\label{eq:parametrized-fractional-knapsack}
\begin{split}
\maximize & \quad  \sum_{k=1}^N  x_{k} w_{\ell} f_{{k}}(f^{\dag}_{{k}}(\alpha/w_{k}))\\
\subject & \quad  \sum_{k=1}^N x_{k}f^{\dag}_{{k}}(\alpha/w_{k}) \le T\\
&\quad x_{k} \in [0,1] ,\quad \forall k \in\until{N}.
\end{split}
\end{align*}

A $2$-factor solution of the knapsack problem with sigmoid utilities is determined in three steps. First, the maximizer $\alpha^*$ of the univariate function $\subscr{F}{LP}$ is determined.  Second, the solution $\boldsymbol x^*(\alpha^*) \in [0,1]^N$ of the $\alpha^*$-parametrized fractional knapsack problem is computed and for each $k\in \until{N}$ for which $x^*_k(\alpha^*)=1$, an amount of resource equal to $f^{\dag}_{{k}}(\alpha/w_{k})$ is allocated, while for every other task no resource is allocated. 
Third, the amount of resource that remains to be allocated in the second step is allocated to the most profitable item with no resource in the second step. 

\end{shaded}
\end{figure}

\setcounter{sideeqn}{\value{equation}}
\setcounter{equation}{\value{eqncnt}}
\renewcommand{\theequation}{\arabic{equation}}

%We denote the receding horizon policy that solves a $N$-horizon certainty-equivalent problem at each stage by $N$-RH policy. 
\section{Design of the Control and Autonomy Management System}
In the previous section, we considered the design of the decision support system. 
%The decision support system design needs the vehicle routing policy $\bs q^\ell$. These inputs are provided by the vehicle routing algorithm. 
We now focus on the other two components of the CAMS, namely, the anomaly detection algorithm and the vehicle routing algorithm.  The anomaly detection algorithm uses the decisions made by the operator to detect an anomalous region and to compute the likelihood of the region being anomalous. This likelihood is utilized by
the vehicle routing algorithm to send the vehicle to an anomalous region with a higher probability. 
%and 
% and vehicle routing algorithm the anomaly detection algorithm, respectively. 
 %In this section, we now focus on the anomaly detection algorithm and vehicle routing algorithm.
\subsection{Anomaly Detection Algorithm}
The decisions made by the human operator about a region being anomalous may be erroneous. The purpose of the anomaly detection algorithm is to reliably decide on a region being anomalous. To this end, we employ statistical quickest change detection algorithms. The quickest change detection algorithm of interest in this paper is the \emph{ensemble CUSUM algorithm}~\cite{VS-FP-FB:11za}.  The ensemble CUSUM algorithm comprises a set of $m$ parallel CUSUM algorithms (one for each region).  We treat the  binary decisions by the operator as Bernoulli random variables  whose distribution is dictated by the performance function of the operator and run the CUSUM algorithm on these decisions to  decide reliably on a region being anomalous. 
%Let the prior belief of the operator on the presence of an anomaly at region $\mc R_{k_\ell}$ just before processing task $\ell$ be $\supscr{\pi}{last}_{k_\ell}$.  
%We treat the binary decisions by the operator on the presence of an anomaly as a Bernoulli random variable and 
%run a 
The standard CUSUM algorithm requires the observations from each region to be independent and identically distributed.  However, the decisions made by the operator do not satisfy these requirements. Therefore, instead of the standard CUSUM algorithm, we resort to the CUSUM like algorithm for dependent observations proposed in~\cite{BC-PW:00}.
We describe the ensemble CUSUM algorithm in Algorithm~\ref{algo:ensemble-CUSUM}.

\begin{algorithm}[ht!]
\begin{algorithmic}[1]
\STATE initialize $\ell:=1$, $\Lambda_k:=0$, for each $k\in \until{m}$
\STATE {\bf if} ${\texttt{dec}}_\ell==1 \text{ and } t_\ell>0$, {\bf  then} 
\[
\Lambda_{k_\ell} := \max \Big\{0, \Lambda_{k_\ell} + \log \frac{f^1_{k_{\ell}}(t_\ell, \pi^{\ell-1}_{k_\ell}) }{1- f^0_{k_{\ell}}(t_\ell,  \pi^{\ell-1}_{k_\ell})} \Big\},
\]
\STATE {\bf else if} ${\texttt{dec}}_\ell==0 \text{ and } t_\ell>0$, {\bf  then} 
\[
\Lambda_{k_\ell} := \max \Big\{0, \Lambda_{k_\ell} + \log \frac{1- f^1_{k_{\ell}}(t_\ell, \pi^{\ell-1}_{k_\ell}) }{f^0_{k_{\ell}}(t_\ell, \pi^{\ell-1}_{k_\ell})} \Big\},
\]

{\it \% detect an anomaly if  a threshold is crossed}
\STATE {\bf if} $\Lambda_{k_\ell} \ge \subscr{\Lambda}{thresh}$, {\bf then}
\STATE \quad declare an anomaly at region $k_\ell$;
\STATE \quad  $\Lambda_{k_\ell}=0$;
\STATE set $\ell=\ell+1$; go to $2:$
\end{algorithmic}
\caption{Ensemble CUSUM Algorithm\label{algo:ensemble-CUSUM}}
\end{algorithm}
%The anomaly detection algorithm treats the decisions by the operator as a Bernoulli random variable whose distribution is dictated by the performance function of the operator. 
The ensemble CUSUM algorithm maintains a statistic $\Lambda_k$ for each region $\mc R_k$, $k\in \until{m}$. The statistic at region $\mc R_k$ is updated using the binary decision of the operator whenever a task from region $\mc R_k$ is processed. 
If the statistic associated with a region crosses a threshold $\subscr{\Lambda}{thresh}$, then the region is declared to be anomalous. The choice of this threshold dictates the accuracy of the detection~\cite{HVP-OH:08}. We assume that once an anomaly has been detected, it is removed, and consequently, the operator's belief about the region being anomalous resets to the default value.

\subsection{Vehicle Routing Algorithm}
We employ a simple vehicle routing algorithm that sends the vehicle to a region with a probability proportional to the likelihood of that region being anomalous. 
In particular, the probability to visit region $\mc R_k$ is initialized to $q_k^0=1/m$ and after processing  each task, the probability to visit region $\mc R_k$ is chosen proportional to $e^{\Lambda_k}/(1+ e^{\Lambda_k})$.  This simple strategy ensures that a region with a high likelihood of being anomalous is visited with a high probability. Moreover, it ensures that each region is visited with a non-zero probability at all times and consequently, an anomalous region is detected in finite time.  

It is noteworthy that such simple vehicle routing algorithm does not take into account the geographic location or the difficulty of detection at each region.   In fact, in the spirit of~\cite{VS-FP-FB:11za}, these factor can be incorporated into the vehicle routing algorithm; however, for simplicity of the presentation, we do not consider these factors here.

\section{A Case Study without Exogenous Factors}
We consider a surveillance mission  involving four regions. The matrix of travel times between the regions is
\[
D= \begin{bmatrix}
0& 22.1422 & 34.4786 &  8.9541 \\
 22.1422 & 0 & 19.3171 & 14.6245 \\
 34.4786 & 19.3171 & 0 &25.5756 \\
  8.9541 & 14.6245 & 25.5756  & 0
\end{bmatrix}. 
\] 
The time to collect information at each region is $10$ units.
The performance of the operator is the same at each region.
Let the drift rate in the DDM associated with the operator be $\mu = \mp 0.3$ for a non-anomalous and an anomalous region, respectively. 
Let the diffusion rate for the DDM associated with the operator be $\sigma=1$. Let the deadline on each task be $40$ units, and let  the importance of each region be the same. Suppose regions $\mc R_1, \mc R_2, \mc R_3$, and $\mc R_4$ become anomalous at time instants $20, 80, 140$, and $200$ units, respectively. 

\begin{figure}[ht!]\footnotesize
\centering
\subfigure[Duration allocation]{
\includegraphics[width=0.8\linewidth]{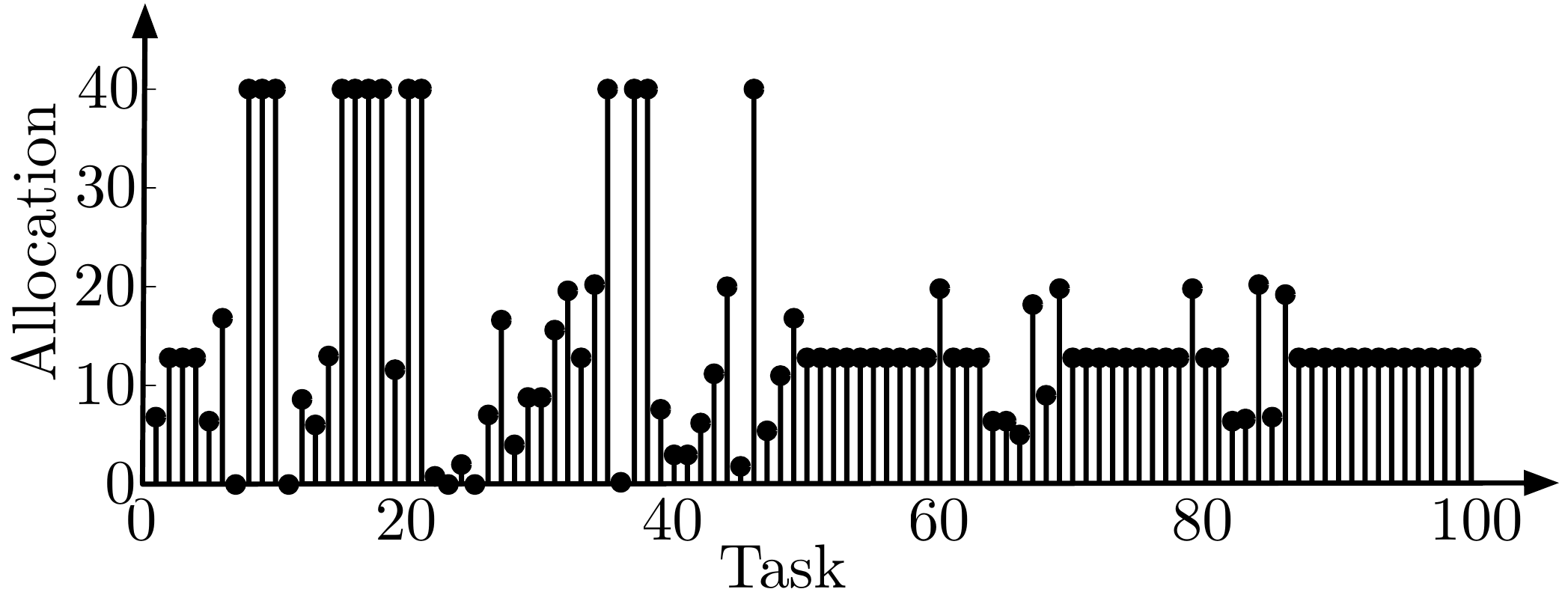} \label{fig:allocation} }\\
\subfigure[Queue length]{
\includegraphics[width=0.8\linewidth]{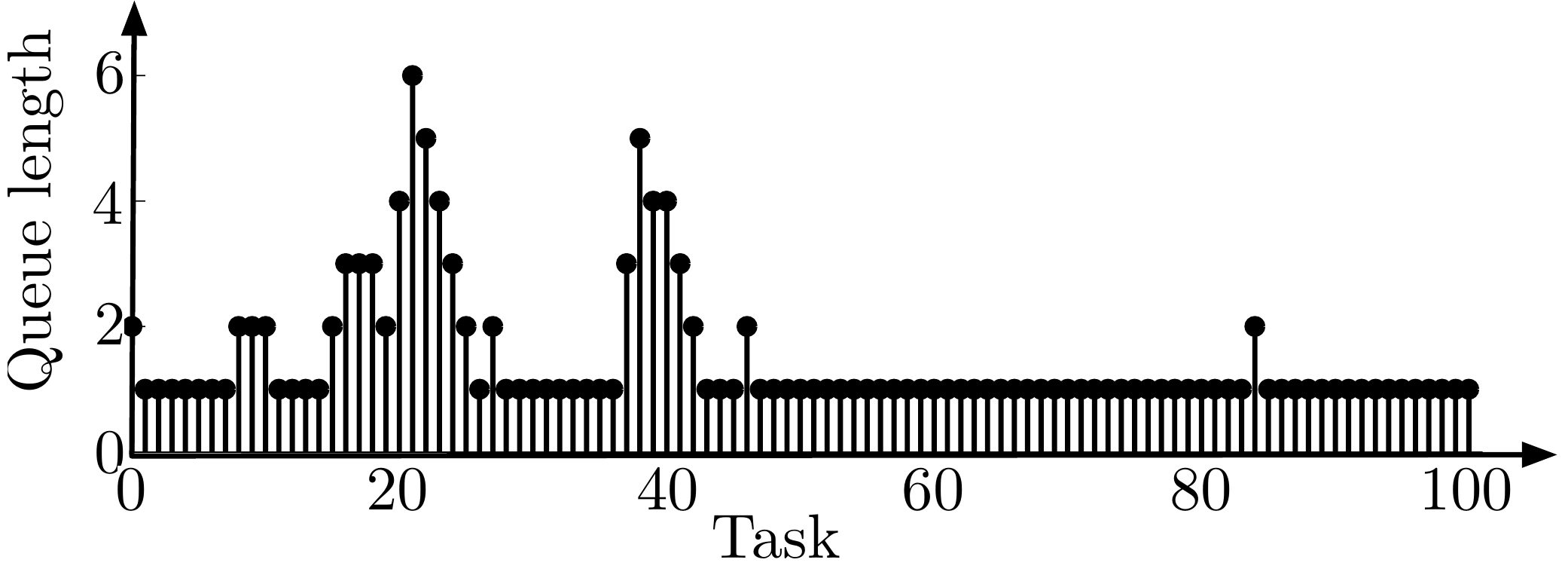}\label{fig:queue}} \\
\subfigure[CUSUM statistics]{
\includegraphics[width=0.8\linewidth]{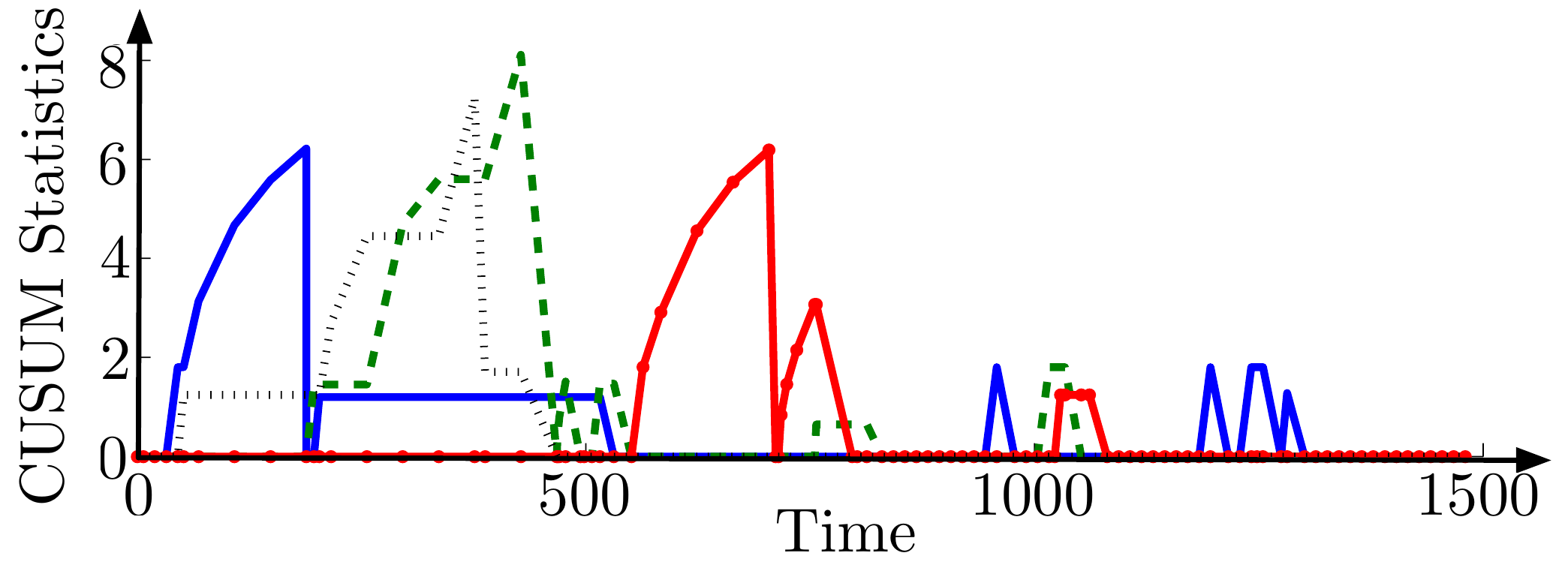}\label{fig:stat}} \\
\subfigure[Region selection probabilities]{
\includegraphics[width=0.8\linewidth]{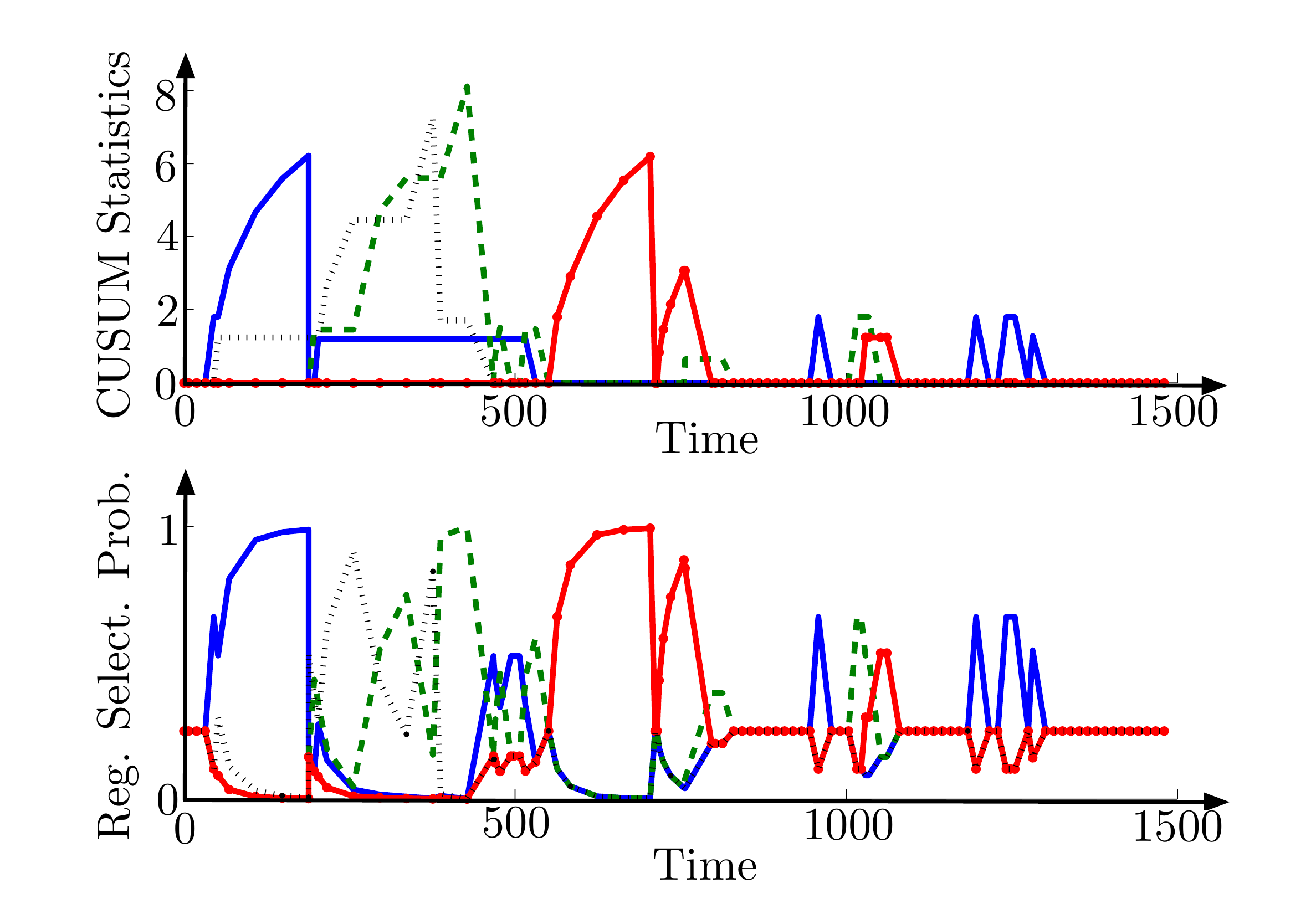}\label{fig:select}} \\
\caption{\footnotesize A sample evolution of the CAMS without exogenous factors for a surveillance mission involving four regions. The solid blue lines, dashed green lines, dotted black lines, and solid red lines with dots represent statistics and region selection probabilities associated with each region. \label{fig:CAMS-sample}}
\end{figure}

The optimization problem~\eqref{eq:maximize-receding-horizon-real-time} with a horizon length $N=5$ is solved before processing each task to determine the optimal allocations for the human operator. A sample evolution of the CAMS is shown in Figure~\ref{fig:CAMS-sample}. As the belief of the operator about the presence of anomaly increases, it becomes more difficult for her to extract more information, and consequently, her reward rate becomes low. In this situation, the allocation policy in Figure~\ref{fig:allocation} allocates a duration equal to the deadline on the task. The queue length in Figure~\ref{fig:queue} increases substantially during this phase and once an anomaly is detected, the allocation policy drops pending tasks in the queue until only one task is present in the queue. The threshold for the CUSUM algorithm is chosen equal to $5$ and once an anomaly is detected the CUSUM statistic in Figure~\ref{fig:stat} resets to zero. 
The region selection policy in Figure~\ref{fig:select} sends the UAV with a high probability to a region with a high likelihood of being anomalous. 

\section{Exogenous Human Factors}
In the previous sections, the human performance model only captured  the effect of evidence aggregation in decision-making.  
There are several exogenous factors, including, situational awareness,  memory retention, fatigue, and boredom that also affect the decision-making performance. In this section, we survey these factor and present a unified model for  human decision-making. The effect of these exogenous factors is typically studied in the free response paradigm for the human operator. Accordingly, we first present the free response paradigm for the DDM.
\subsubsection{Drift-diffusion model and the free response paradigm}
\begin{figure*}[ht!]
\centering  \scriptsize
\subfigure[Free Response Paradigm]{
\includegraphics[width=0.3\textwidth]{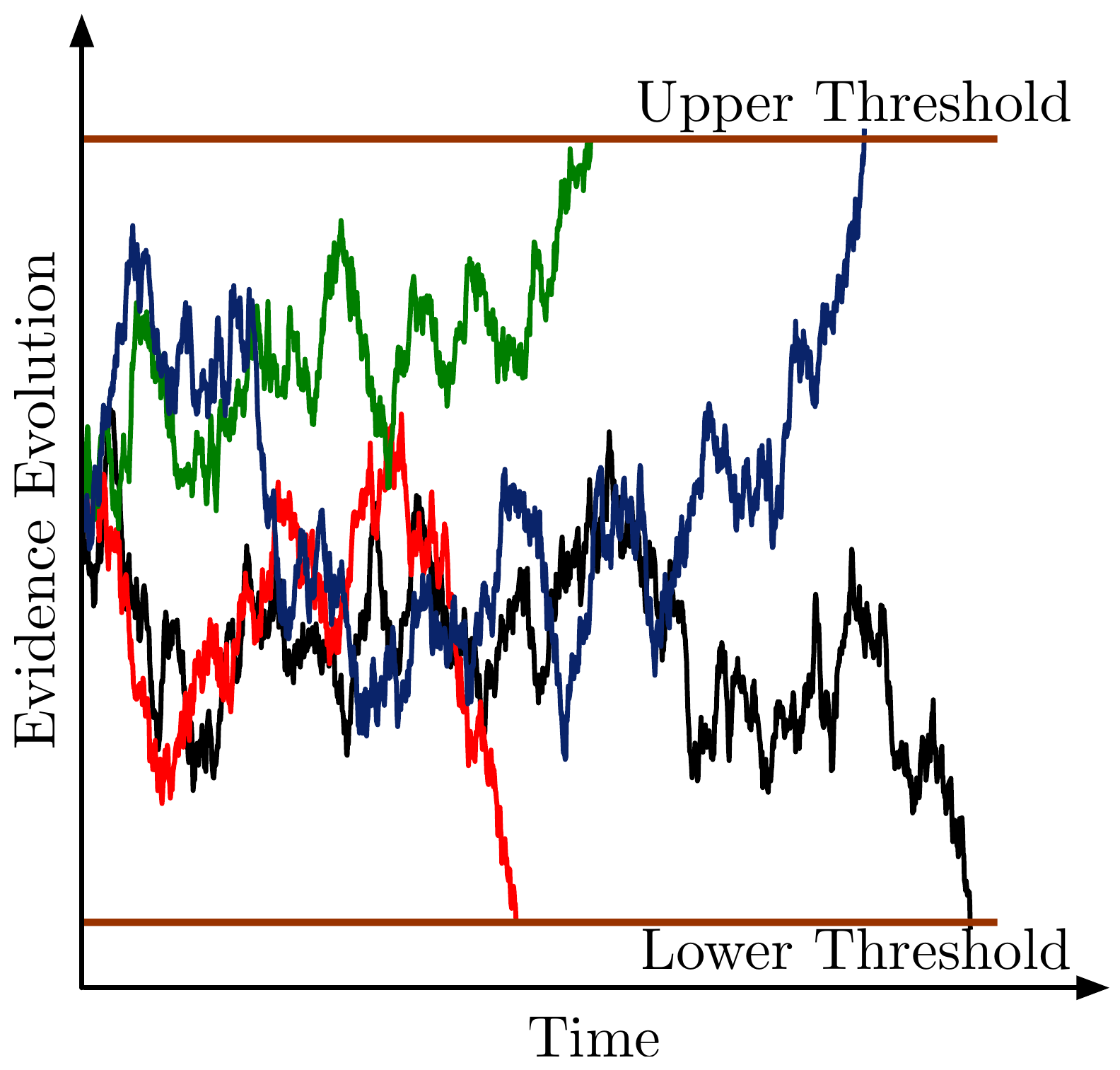} \label{fig:free-response-ddm}
}  
\subfigure[Empirical Reaction Times]{
\includegraphics[width=0.3\textwidth]{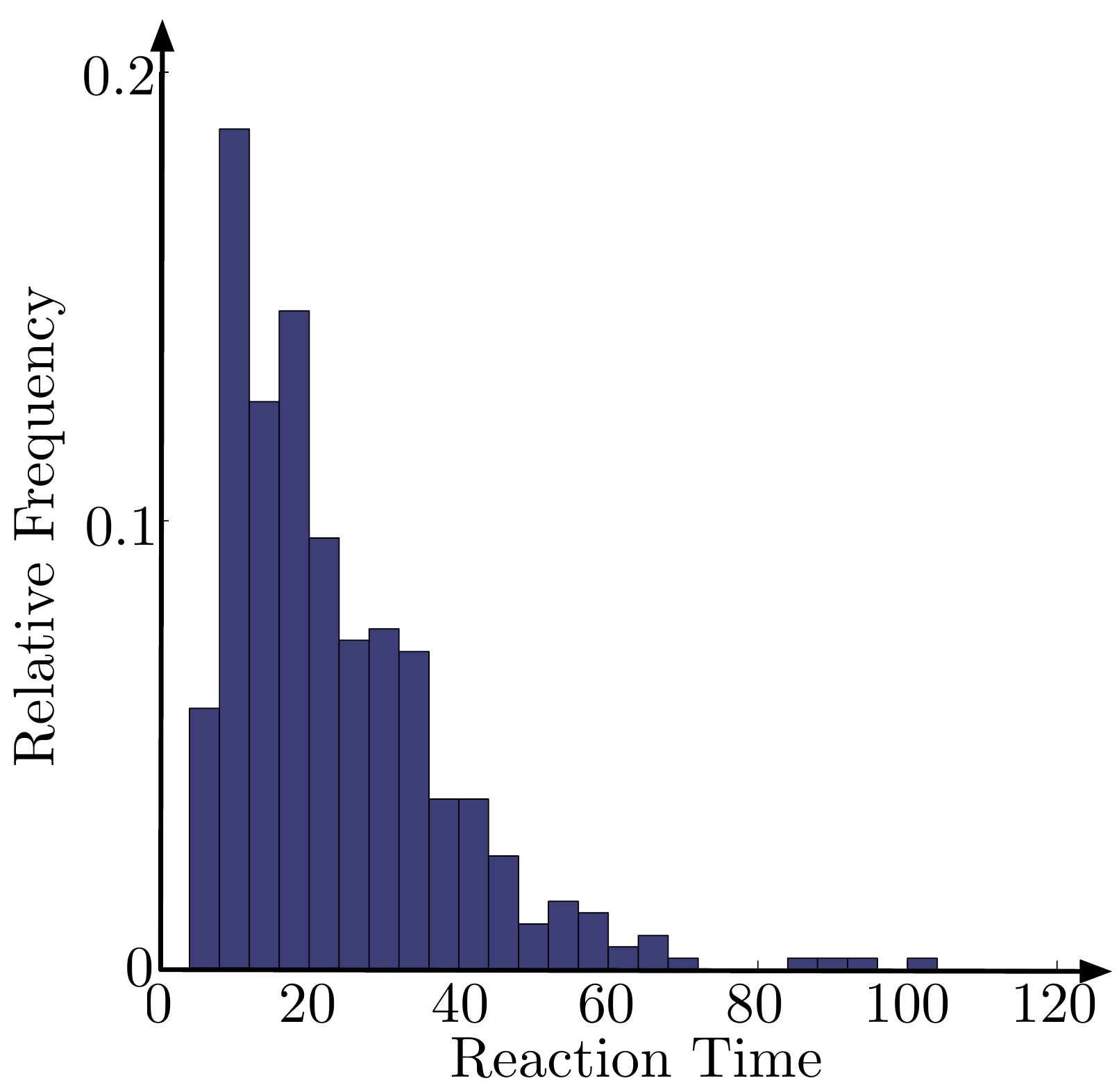} \label{fig:RT_data}
} 
\subfigure[Reaction Times under the Free Response Paradigm]{
\includegraphics[width=0.3\textwidth]{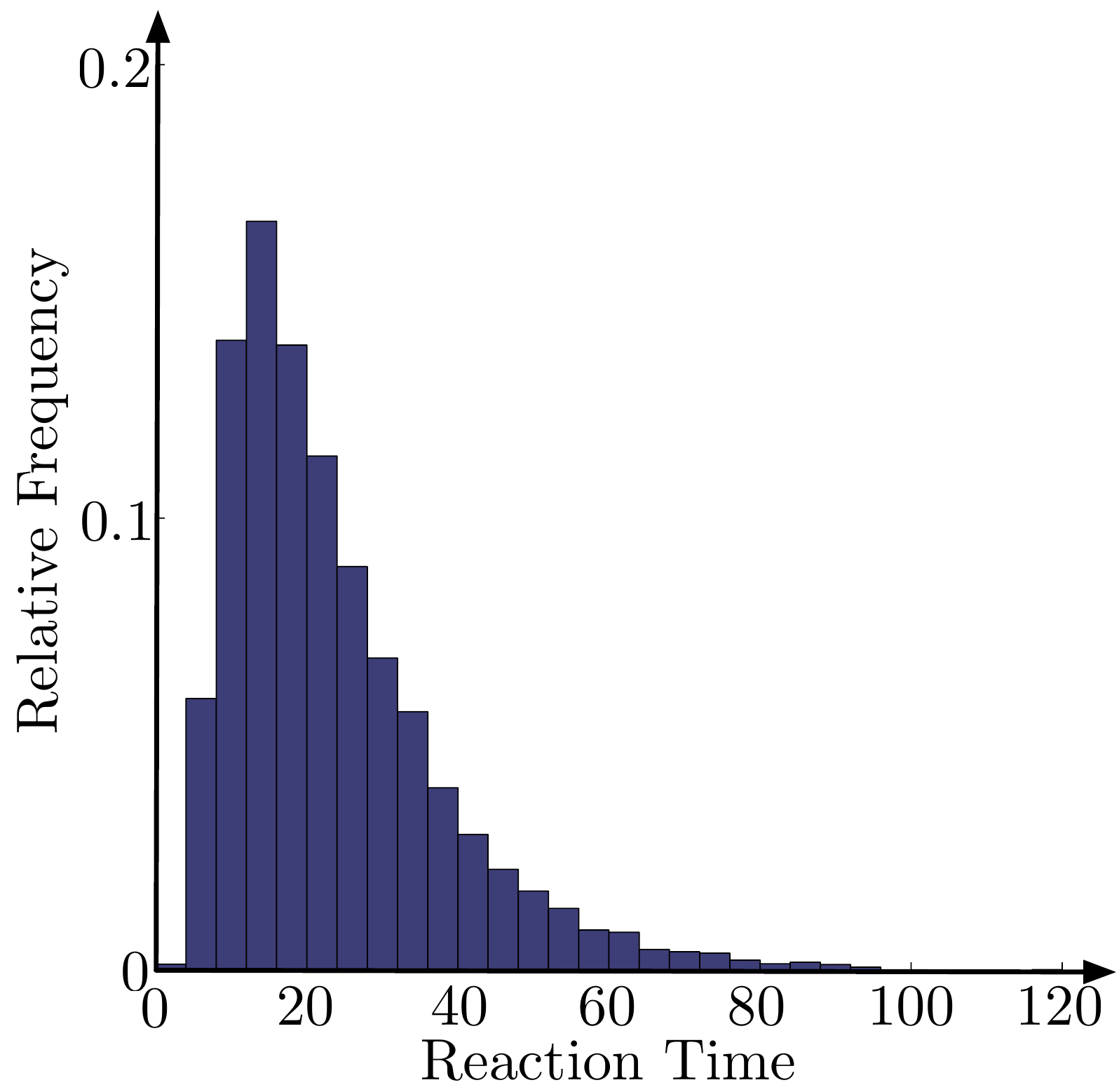} \label{fig:RT_ddm}
}
\subfigure[Waiting time due to lack of situational awareness]{
\includegraphics[width=0.3\textwidth]{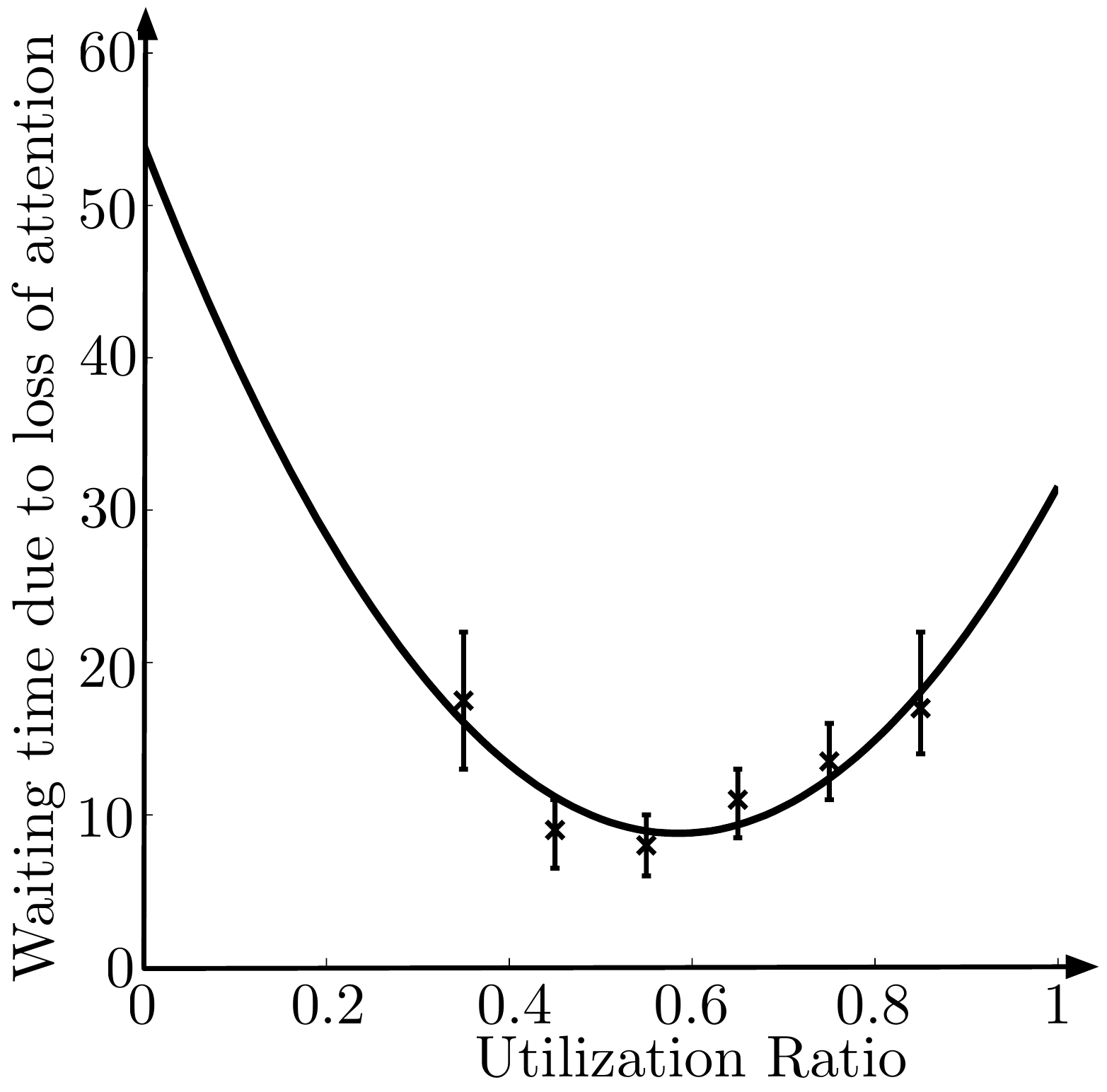} \label{fig:wtsa}
}
\subfigure[Effect of Fatigue]{
\includegraphics[width=0.3\textwidth]{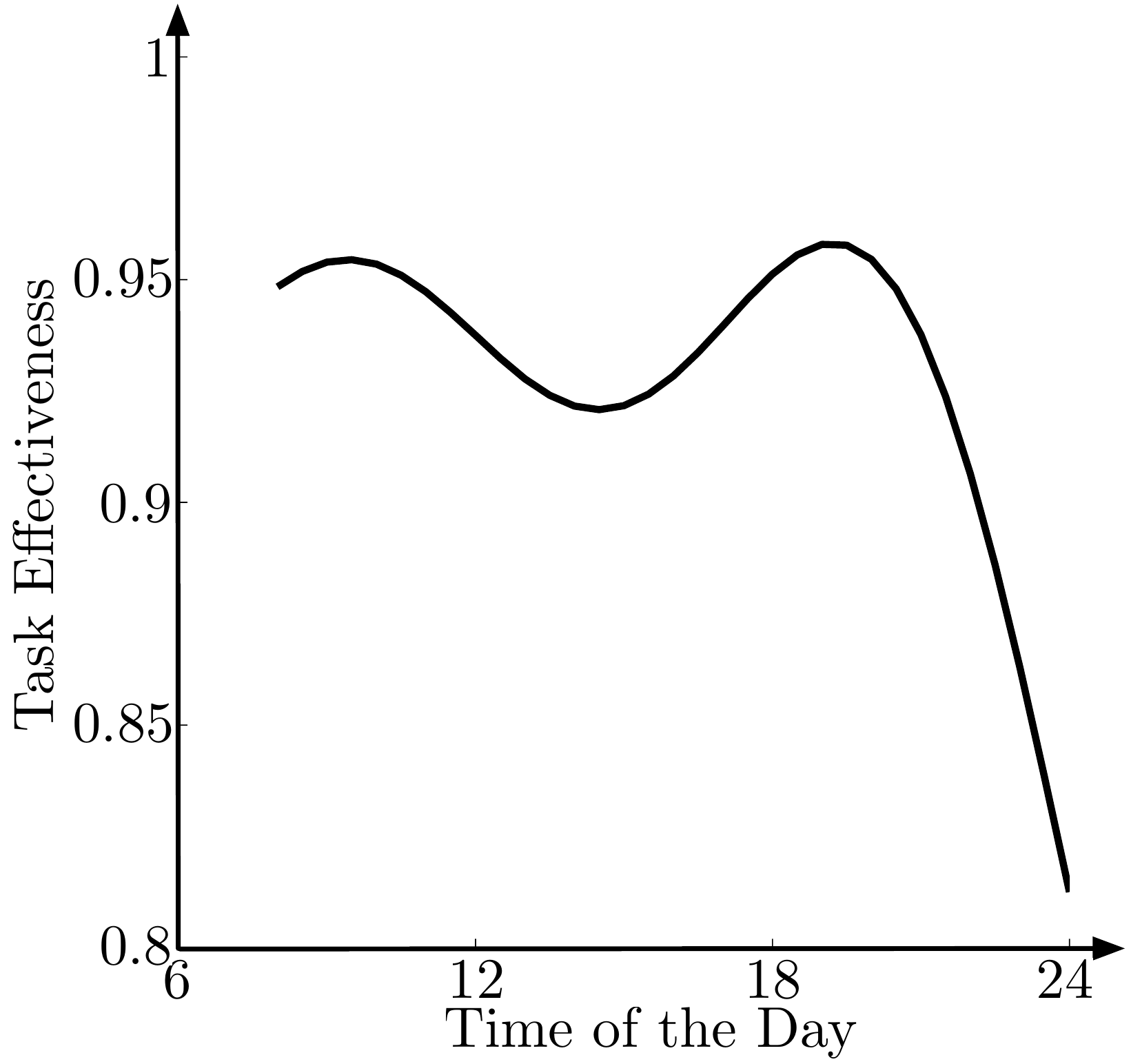} \label{fig:safte}
} 
\subfigure[Memory Retention]{
\includegraphics[width=0.3\textwidth]{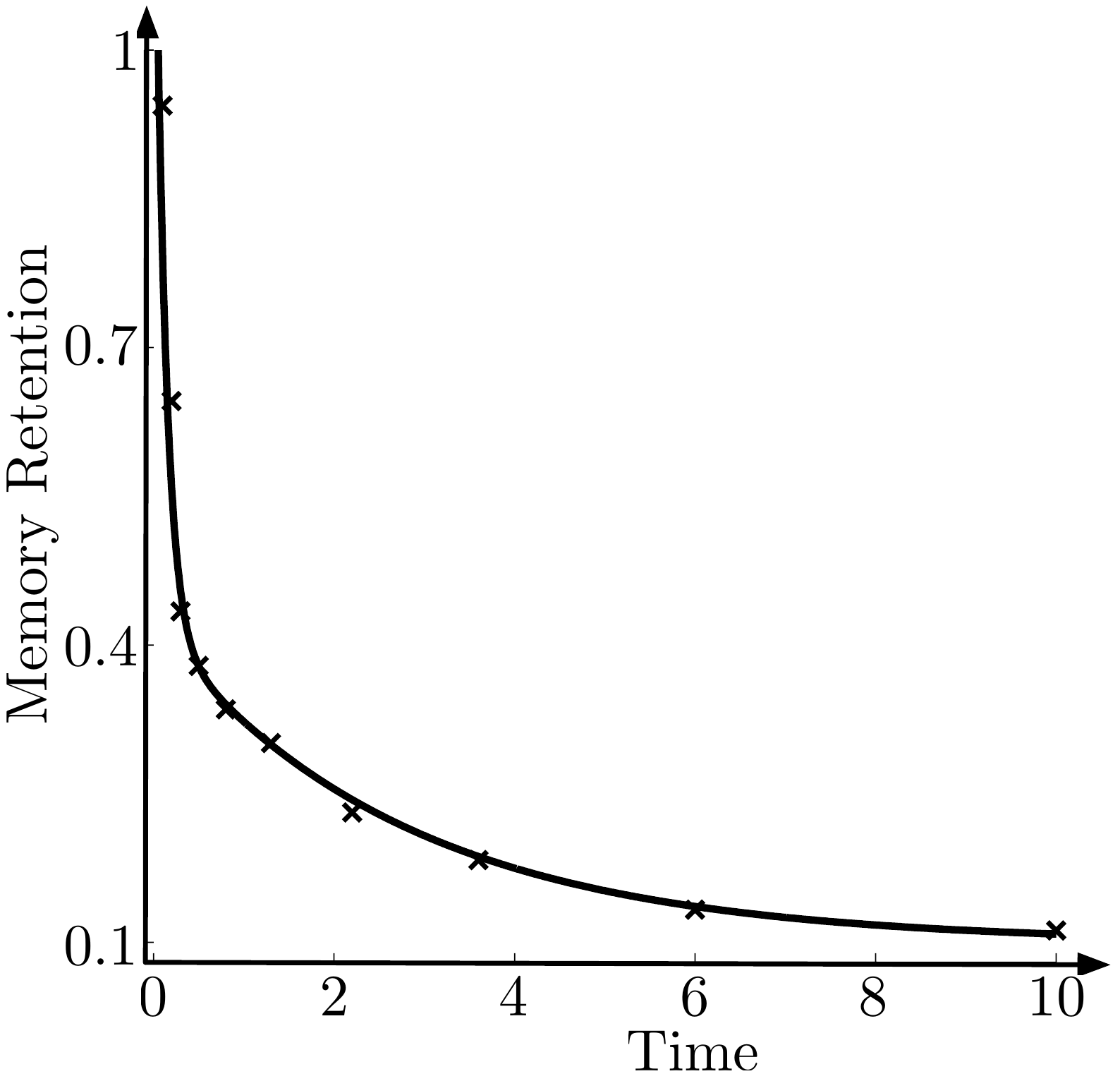}\label{fig:retention}
}
\caption{
(a) The free response paradigm for decision-making. The evidence evolves according to the DDM and the operator makes a decision whenever one of the thresholds is crossed. (b) Empirical reaction time data taken from~\cite{KS-CN-TT-EF:08}. (c) Decision time distribution for a DDM under the free response paradigm. The parameters were estimated using the data in~\cite{KS-CN-TT-EF:08}.  (d) Delay in processing a task due to the lack of situational awareness. (e) Task effectiveness of a operator who wakes up at $6$am after $6$ hours of sleep. (f) Memory retention as a function of time.}
\end{figure*}

In the free response paradigm, the operator takes her own time to decide on an alternative. For the DDM~\eqref{eq:ddm}, 
the free response paradigm is modeled via two thresholds (positive and negative) and the operator decides in favor of the first/second alternative if the positive (negative) threshold is crossed from below (above). 
A typical evolution of the DDM under free response paradigm is shown in Figure~\ref{fig:free-response-ddm}.
If the two alternatives are equally likely, then the two thresholds are chosen symmetrically.  Let $\pm \eta \in \real$ be symmetrically chosen thresholds. The expected decision time ($\subscr{T}{decision}$) under the free response paradigm is
\[
\subscr{T}{decision}= \frac{\eta}{\mu} \tanh \frac{\mu \eta}{\sigma^2} + \frac{2 \eta (1 - e^{-2x_0 \mu/\sigma^2})}{ \mu(e^{2\eta \mu /\sigma^2}-e^{-2\eta \mu /\sigma^2}) } -\frac{x_0}{\mu}.
\]
The reaction time on a task is  $\subscr{T}{decision}+\subscr{T}{motor}$, where $\subscr{T}{motor} \in \real_{>0}$ is the time taken by sensory and motor processes unrelated to the decision process. For simplicity, we treat $\subscr{T}{motor}$ as a deterministic quantity in this paper.

The choice of the threshold is dictated by the speed-accuracy trade-off. 
There are two particular criteria to capture the speed accuracy
trade-off: (i)  the Bayes risk and (ii) the reward rate~\cite{RB-EB-etal:06}. We
focus on the Bayes risk criterion in this paper. The Bayes risk is defined
by $\texttt{BR}= \xi_1 \subscr{T}{decision}+ \xi_2 \subscr{\prob}{error}$, where $\xi_1, \xi_2$
are the cost per unit delay in decision and the cost of error, respectively,
and $\subscr{\prob}{error}$ is the error rate. For the DDM, the minimization of
the Bayes risk yields the following transcendental equation for the
threshold~\cite{RB-EB-etal:06}:
\[
\frac{\xi_2}{\xi_1} \frac{2 \mu^2}{\sigma^2} -\frac{4 \mu \eta}{\sigma^2} +e^{-(2\mu \eta/\sigma^2)}-e^{(2\mu \eta/\sigma^2)}=0.
\]
This transcendental equation can be solved numerically but for the purpose
of this paper, we treat the limiting solution of this equation as the true
solution. In particular, if $\mu$ is low or $\sigma$ is high, then the
threshold is $\eta = \mu \xi_2/4\xi_1$.  It is noteworthy that $\xi_1$ and $\xi_2$
can be estimated from the empirical data (see~\cite{RB-EB-etal:06} and
references therein).

\subsubsection{Yerkes-Dodson law and situational awareness}

The Yerkes-Dodson law~\cite{RMY-JDD:08, CDW-JGH:00} is a classical model that captures the performance of an operator as a unimodal function of their arousal level. The arousal level of an operator is a measure of her stress.
The Yerkes-Dodson law
suggests that there exists an optimal level of arousal at which the operator's performance is the best. Moreover, the optimal arousal level decreases with the difficulty of the task.

Situational awareness of an operator captures the level of her spatial and temporal perception of the environment. The Yerkes-Dodson law is used to model the situation awareness of  the operator as a function of her workload~\cite{KS-CN-TT-EF:08, KS-EF:10b, CN-BM-JWC-MLC:08}. In particular, the 
workload of the operator is chosen as a measure of the arousal level, and the situational awareness
is modeled as a unimodal function of the workload. 

The reaction time of the operator decreases with increasing situational awareness and accordingly, the expected reaction time is modeled as a convex function of the workload~\cite{KS-EF:10b}. It has been argued in~\cite{MC-CN-etal:07} that the lack of situational awareness results in a larger waiting time for the tasks, i.e., the operator takes more time to start working on a task. This suggests that the lack of situational awareness does not affect the decision-making process and only increases the sensory and motor time $\subscr{T}{motor}$.  Accordingly, we model the sensory and motor time as a convex function of the workload and treat the expected decision time as a constant function of the workload. 

The workload is modeled as the utilization ratio (the fraction of recent history during which the operator was busy) and the utilization ratio $u$ is captured through the following differential equation
\[
\dot u(t) = \frac{b(t)-u(t)}{\tau}, \quad u(0)=u_0,
\]
where $\map{b}{\real_{>0}}{\{0,1\}}$ is the binary function of time $t$ that represents if the operator is idle or busy, $\tau\in \real_{>0}$ is the sensitivity of the operator, and $u_0 \in [0,1]$ is the initial utilization ratio of the operator~\cite{KS-EF:10b}. 

We denote the sensory and motor time as a function of the utilization ratio by $\map{\subscr{T}{motor}}{[0,1]}{\real_{\ge 0}}$. The expected reaction time ($\subscr{T}{decision}+\subscr{T}{motor}$) is also a convex function of the utilization ratio and such a relation  has been empirically validated in~\cite{KS-EF:10b}. Moreover, it has been noted in~\cite{CN-BM-JWC-MLC:08} that these utilization ratio based models only capture the effect of the workload and do not capture training and fatigue effects.

\subsubsection{Fatigue, sleep cycle and SAFTE model}
Fatigue is the feeling of tiredness and bodily discomfort after prolonged activity~\cite{GM-DRD-etal:00}. Fatigue is known to have detrimental effects on operator's performance. Several models for fatigue have been proposed; a summary is presented in~\cite{MMM-SM-etal:04}. In this paper, we adopt Sleep  Activity Fatigue Task Efficiency (SAFTE) model proposed in~\cite{SRH-DPR-etal:04}. This model considers three important processes, namely, homeostatic regulation of wakefulness, circadian process, and sleep inertia to model cognitive performance as a function of sleep deprivation.
The SAFTE model assumes that a fully rested operator has a finite maximum cognitive capacity called the reservoir capacity $R_c$. While the operator is awake, the cognitive reservoir is depleted linearly with time. The reservoir replenishes when the operator sleeps and the replenishment depends on the circadian process (time of the day) and the current reservoir level.
We consider the same form of the SAFTE model as in~\cite{NDP-KAM:12}. The SAFTE model determines the task effectiveness as
\begin{multline*}
\texttt{TE}= 100 \frac{R_c-60K T_a}{R_c} + \Big( a_1 +a_2 \frac{60 K T_a}{R_c}\Big) \\
 \Big[ \cos \Big( \frac{2\pi}{24}(T_d-p) \Big)+\beta \cos \Big(\frac{4\pi}{24}(T_d-p-p')\Big)\Big],
\end{multline*}
where $T_a$ is the number of hours the operator has been awake, $T_d$ is the time of the day in hours, $K$ is reservoir drain rate due to wakefulness, $a_1, a_2, \beta \in \real$ are constants, $p$ is the time of the peak in the $24$h circadian rhythm and $p'$ is the relative time of the $12$h peak. The default values of the parameters in~\cite{SRH-DPR-etal:04} are $R_c= 2880$ units, $K= 0.5$ units per minute, $a_1=7$, $a_2=5$, $\beta=0.5$, $p=18$ hours, and $p'=3$ hours.
The task effectiveness is a measure of the efficiency of the operator. In particular, if the reaction time of a fully rested operator is \texttt{RT}, then the reaction time of the fatigued operator is $\texttt{RT}/ \texttt{TE}$.

\subsubsection{Forgetting/Retention curve}
The forgetting/retention curve determines the fraction of newly acquired information the operator remembers over time. 
The modeling of the forgetting curve has been a debated topic. Traditionally, the forgetting curve has been modeled as an exponential decay~\cite{HE:13}.  Anderson~\etal~\cite{JRA-LJS:91} argue that the forgetting curve should be modeled by a power law function.  Rubin~\etal~\cite{DCR-SH-AW:99} model the forgetting curve as a sum of two exponential functions and a constant function. They claim that one of the exponential terms describes the working memory, and the two remaining terms describe the long term memory.  For the purpose of this paper, the functional form of the forgetting curve is not important. We assume that the fraction of memory retained is known as a function of time and we denote it by $\map{\texttt{rem}}{\real_{\ge 0}}{[0,1]}$. 

\subsubsection{Motivation and boredom curve}
Boredom in an important factor in human-in-the-loop systems that is mostly neglected in design~\cite{CDF:93} and leads to operator's lack of interest in their activity.  To the best of our knowledge, there has not been a significant work in the mathematical modeling of  boredom. A heuristic model for the decay in performance due to boredom has been developed in~\cite{NA-SZ-ML:10}. The authors argue that the lack of motivation can be thought of as boredom, and hypothesize an exponential decay model based on the classical motivation theory. Although the analytic approach in this paper can easily incorporate boredom models, due to lack of empirical evidence validating boredom models, we do not consider them in the remainder of the paper.

%\begin{figure*}
%\centering
%  \hfill
%
%}
%\caption{}
%\end{figure*}

\subsection{A Unified Operator Model}

We now develop a unified operator model that blends the aforementioned models and determines the performance of the operator as a function of the duration allocation, her utilization ratio,  her sleep schedule, and her forgetfulness. 
Without loss of generality, we assume that the parameters for the decision dynamics of the operator were determined when the operator was fully rested.
Let the estimated drift and diffusion rates for the DDM associated with the operator be $\mu$ and $\sigma$, respectively. It follows from the definition of the task effectiveness that the expected decision time of the fatigued operator is
\[
\frac{1}{\texttt{TE}} \frac{\eta}{\mu} \tanh\frac{\mu \eta}{\sigma^2},
\]
where $\texttt{TE}$ is the task effectiveness obtained using the SAFTE model. 
We assume that the operator fatigue  affects their drift rate while the diffusion rate remains the same. Let $\supscr{\mu}{eff}$ be the drift rate of the DDM associated with the fatigued operator and let $\supscr{\eta}{eff}$ be the associated threshold. 
It is known that under the limiting approximation $\supscr{\eta}{eff} = \supscr{\mu}{eff} \xi_2/4\xi_1$. It follows that
\begin{equation}\label{eq:modified-drift}
 \tanh \frac{{\supscr{\mu}{eff}}^2 \xi_2}{4 \xi_1 \sigma^2} =\frac{1}{\texttt{TE}} \tanh\frac{\mu \eta}{\sigma^2}.
\end{equation}
Equation~\eqref{eq:modified-drift} yields the following expression for the drift rate of a fatigued operator:
\[
\supscr{\mu}{eff}= \sqrt{ \frac{2 \xi_1 \sigma^2}{\xi_2} \log \Big(  \frac{\texttt{TE}+ \tanh(\mu^2 \xi_2/ 4 \xi_1 \sigma^2)}{\texttt{TE}-  \tanh (\mu^2 \xi_2/ 4 \xi_1 \sigma^2)}\Big)}.
\]
Similarly, the sensory and motor time $\subscr{T}{motor}$ for a fatigued operator would be $\subscr{T}{motor}/\texttt{TE}$. 

We now consider the retention model. Consider the current decision-making task that the operator processes.  Suppose a decision-making task from the same class as the current task was processed $\subscr{T}{last}$ time earlier and let the evidence after that processing was $\subscr{\pi}{last}$. Thus, according to the forgetting model, the initial condition for the DDM  associated with the current task is  $\subscr{x}{init} = (\sigma^2 \log (\subscr{\pi}{last}/(1-\subscr{\pi}{last}))/2\supscr{\mu}{eff} ) \texttt{rem}(\subscr{T}{last})$. 
%$\subscr{\pi}{last} $. This maps to an initial state 

The decision support system suggests the operator to allocate a given duration to the current task. Such a situation corresponds to the interrogation paradigm of decision-making.  A unified performance function for the operator can be determined by using the effective drift rate, the effective motor time, and the effective initial condition for the associated DDM. 
Thus, the performance function $\map{f^1}{\real_{\ge 0} \times [0,1]\times [0,1] \times \real_{\ge 0} \times [0,1]}{[0,1]}$ of the operator on a task from an anomalous region is
\begin{equation}\label{eq:performance-function}
f^1(t, u, \texttt{TE}, \subscr{T}{last}, \subscr{\pi}{last}) = 1- \Phi\Big( \frac{\nu -\subscr{\mu}{eff}(t-\supscr{T}{wait})^+ -\subscr{x}{init}}{\sigma \sqrt{(t-\supscr{T}{wait})^+}}\Big),
\end{equation}
where $\supscr{T}{wait}= \subscr{T}{motor}(u)/\texttt{TE}$ and $(\cdot)^+ := \max\{0, \cdot\}$. The performance function of the operator on  a task from a non-anomalous region can be defined similarly. The performance function in equation~\eqref{eq:performance-function} assumes that no evidence is collected during the sensory and motor time. 

In summary, the above performance function can be interpreted in the following way. The belief of the operator about a region being anomalous renders an initial condition to the associated DDM. This belief is acquired while the operator processed the task from the same region last time. Moreover, over the course of time the operator forgets the acquired belief. Accordingly, the belief acquired after  a task from the current region was last processed is discounted using the forgetting model, and the discounted belief is used as the  
the initial condition for the DDM associated with the current task. 
Furthermore, the operator fatigue is captured by modifying the drift rates of the DDMs associated with different tasks  to match the first moment of the reaction times suggested by the SAFTE model.

\subsection{Design of the CAMS using unified operator model}

We now design the CAMS using the unified operator
model~\eqref{eq:performance-function}.  The design proceeds analogously to the
design with the simplified operator~\eqref{eq:performance-func}. In this
section we focus only on the design of the decision support system.  We
study the problem using the certainty-equivalent receding-horizon
framework.  We assume that the parameters of the operator remain constant while a
given task is processed.  The first exogenous factor we consider is the
utilization ratio.  The utilization ratio of the operator affects their
performance and can be controlled by introducing an idling time for the
operator. The idling time reduces the fraction of time the operator is busy
and hence reduces the utilization ratio of the operator. Let the
utilization ratio of the operator be $u_\ell$ before processing the
$\ell$-th task.  If the operator processes the $\ell$-th task for a duration $t_\ell$
and remains idle for a duration $\delta_\ell$ after processing the task,
then the utilization ratio after processing the $\ell$-th task  is
\[
u_{\ell+1} =  \Big(1-e^{-{t_\ell}/ {\tau}} + u_{\ell} e^{-{t_\ell}/ {\tau}}\Big) e^{-{\delta_\ell}/ {\tau}}. 
\]

With the incorporation of the utilization ratio dynamics, the problem has two state variables, namely, the queue length and the utilization ratio. There are two control variables as well, namely, the duration allocation to a task, and the rest time that follows the processing of the task. The duration allocation controls the decision-making performance and the rest time controls the utilization ratio of the operator. 
To capture the memory effects, the optimization problem also needs to keep track of the time when a task
from a region was last processed and the time when a task from a region will be next processed.
Thus, the number of state variables increases further. 
The increased dimension of the state space and action space would make the computation of the solution to the finite-horizon problem associated with the certainty-equivalent receding-horizon control intractable. However, for this finite-horizon optimization problem some of the variables can be assumed to remain constant, and this reduces the computational complexity.

The utilization ratio affects the sensory and motor time of the operator. The sensory and motor time is high only at extreme values of the utilization ratio. Therefore, if we control the utilization ratio of the operator such that it always lies within an interval around its optimal value, then it would result in efficient sensory and motor times for the operator. In particular, a threshold $\subscr{u}{th}$ can be chosen such that the sensory and motor time does not vary much in the interval $[\subscr{u}{opt},\subscr{u}{th}]$, where $\subscr{u}{opt}$ is the utilization ratio at which the minimum sensory and motor time is achieved.
Moreover, after processing each the $\ell$-th task, if the utilization ratio $\bar u_{\ell}$ is above a threshold value $\subscr{u}{th}$, then a idling time $\delta_\ell = \tau \log (\bar u_{\ell} /\subscr{u}{opt})$ is suggested to the operator such that the utilization ratio of the operator goes down to its optimal value $\subscr{u}{opt}$. 
Such a threshold based control of the utilization ratio ensures that we can treat the motor time as a constant over the finite prediction horizon and consequently, we do not need to consider the idling time as a decision variable in the finite-horizon problem. 

The fatigue and sleep cycle models are slow timescale models, i.e., the performance significantly changes  only after an elapse of hours of duration, while the evidence aggregation process in the decision-making performance is of the order of seconds.
Accordingly, for the finite-horizon problem associated with the receding-horizon policy, we treat the task effectiveness \texttt{TE} as a constant.  Similarly, the memory models are of the order of minutes and accordingly, the operator's belief can be assumed to be a constant over the finite prediction horizon. 
Thus, the fatigue and the memory effects need not be considered in the finite-horizon problem associated with the certainty-equivalent receding-horizon control. 

%, we do not need to keep track of the time when the task from a region was last processed. 
%the tracking of time when the task would be processed is not needed. 
%, while  and the fatigue models are of order of hours. Therefore, over the prediction horizon operator's belief and task effectiveness can be assumed to be constant. 
%This makes the finite horizon problem associated with certainty-equivalent receding horizon control a dynamic program with univariate state and control variables. 
%

Let $ \supscr{T}{last}_{k}$ be the time since the task from region $\mc R_{k}$ was last processed and let $\supscr{\pi}{last}_{k}$ be the belief of the operator about region $\mc R_{k}$ being anomalous after  a task from region $\mc R_k$  was last processed.
%Let $\supscr{T}{visit}_{k}$ be the last time when operator processed an evidence from region $\mc R_k$. 
%Let $\supscr{T}{glob}$ be the current time. It follows that $\supscr{T}{last}_{k} =\supscr{T}{glob}- \supscr{T}{visit}_{k}$.
 In view of the above discussion, for the finite-horizon problem to be
 solved when the $\ell$-th task  in the queue is processed, the performance function
 $\map{f_{k}^{\ell}}{\real_{\ge 0}}{\real}$ associated with a task from
 region $\mc R_k$ is defined by
\begin{multline}\label{eq:certain-performance}
f_{k}^\ell (t) =
 (1- \pi^{\ell-1}_{k}) f_{k}^0 (t, u_{\ell-1}, \texttt{TE}_\ell, \supscr{T}{last}_k, \supscr{\pi}{last}_k)  \\
+ \pi^{\ell-1}_{k} f_{k}^1 (t, u_{\ell-1}, \texttt{TE}_\ell, \supscr{T}{last}_k, \supscr{\pi_k}{last}),
\end{multline}
where $\texttt{TE}_\ell$ is the task effectiveness while processing the $\ell$-th task, and $f_{k}^0$ and $f_{k}^1$ are computed as in
equation~\eqref{eq:performance-function}. After processing the $\ell$-th task, the
utilization ratio, and the task effectiveness are updated using associated models,
while the belief at region $\mc R_j$ is updated using Bayes rule as
\[
\bar \pi_{j}^{\ell} =
\begin{cases}
 \frac{\supscr{\pi_{j}}{rem} \prob({\texttt{dec}}_\ell| H^1_{k})}{(1- \supscr{\pi_{j}}{rem}) \prob({\texttt{dec}}_\ell| H^0_{k_\ell})+ \supscr{\pi_{j}}{rem} \prob({\texttt{dec}}_\ell| H^1_{k_\ell})}, & \text{if } j=k_{\ell}, \\
 \supscr{\pi_{j}}{rem}, & \text{otherwise,}
 \end{cases}
\]
where $H^0_k$ and $H^1_k$ denote the hypothesis that region $\mc R_k$ is non-anomalous and anomalous, respectively, and ${\texttt{dec}}_\ell \in \{0,1\}$ is the operator's decision, $\supscr{\pi_{k}}{rem}$ is the operator's belief about  region $\mc R_k$ being anomalous after accounting for the memory retention effects and is defined by
\[
\supscr{\pi_{k}}{rem} = \frac{ \exp \big \{\log\big( \frac{\supscr{\pi_{k}}{last}}{1- \supscr{\pi_{k}}{last}}\big) \texttt{rem}(\supscr{T}{last}_{k} )\big \}}
{1+ \exp \big \{\log\big( \frac{\supscr{\pi_{k}}{last}}{1- \supscr{\pi_{k}}{last}}\big) \texttt{rem}(\supscr{T}{last}_{k} )\big \}}, 
\]
and $\prob({\texttt{dec}}_\ell|\cdot)$ is determined from the performance function of the operator, e.g., 
\[
\prob({\texttt{dec}}_\ell=1| H^1_{k_\ell}) = f^1_{k_\ell} (t_\ell, u_{\ell-1}, \texttt{TE}, \supscr{\pi_{k_\ell}}{last}, \supscr{T_{k_\ell}}{last}).
\]
Similar to the simplified operator model case, the operator resets her belief to a threshold value if her belief is smaller than the threshold. Accordingly, her belief about region $\mc R_k$ being anomalous after processing the $\ell$-th task is
\[
\pi_{j}^{\ell} = \max \{0.5, \bar \pi_{j}^{\ell}\}.
\]
Further, after processing the $\ell$-th task,  $\supscr{\pi}{last}_{k_\ell}$ is updated to $\pi_{k_\ell}^{\ell} $.

The receding-horizon duration allocation policy now follows similarly to the simplified operator model case. The duration allocation to the $\ell$-th task  is determined by solving a problem of the form~\eqref{eq:maximize-receding-horizon-real-time} in which the performance function for each task is of the form~\eqref{eq:certain-performance}.
Again, the finite-horizon problem is a dynamic program with univariate state and control variables and an efficient solution can be determined using the backward induction algorithm on the discretized state and action space.

The design of the anomaly detection algorithm and the vehicle routing policy follow similarly to the simplified operator model case. 

\section{A Case Study with Exogenous Factors}
We now present a case study on the design of the CAMS with exogenous
factors. We choose the same parameters and the same time instances for the appearance of
anomalies as in
the case without exogenous factors. Additionally, the deadline on each task
is $60$ units, the operator sensitivity is $100$, and the motor time is
$\supscr{T}{motor}(u)=54 -155u + 132u^2 -9$. We choose the retention
function of the form in~\cite{DCR-SH-AW:99} given by
$\texttt{rem}(t)=\min(1, a_1 \exp(-10t/1.15)+ a_2 \exp(-10t/27.55)+a_3)$,
where $a_1= 4.6$, $a_2=1.5$, and $a_3=0.1$.

\begin{figure*}[ht!]\footnotesize
\centering
\subfigure[Duration allocation]{
\includegraphics[width=0.4\linewidth]{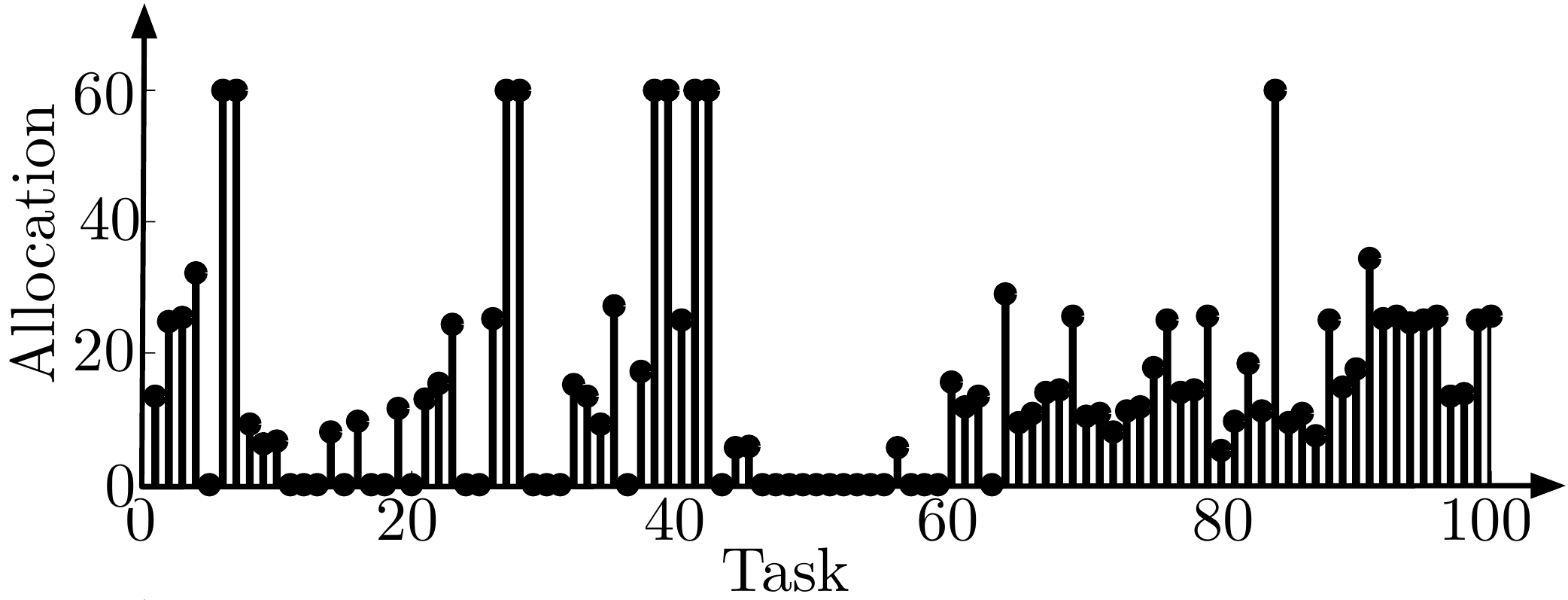} \label{fig:allocation-full} }\qquad
\subfigure[Rest Time]{
\includegraphics[width=0.4\linewidth]{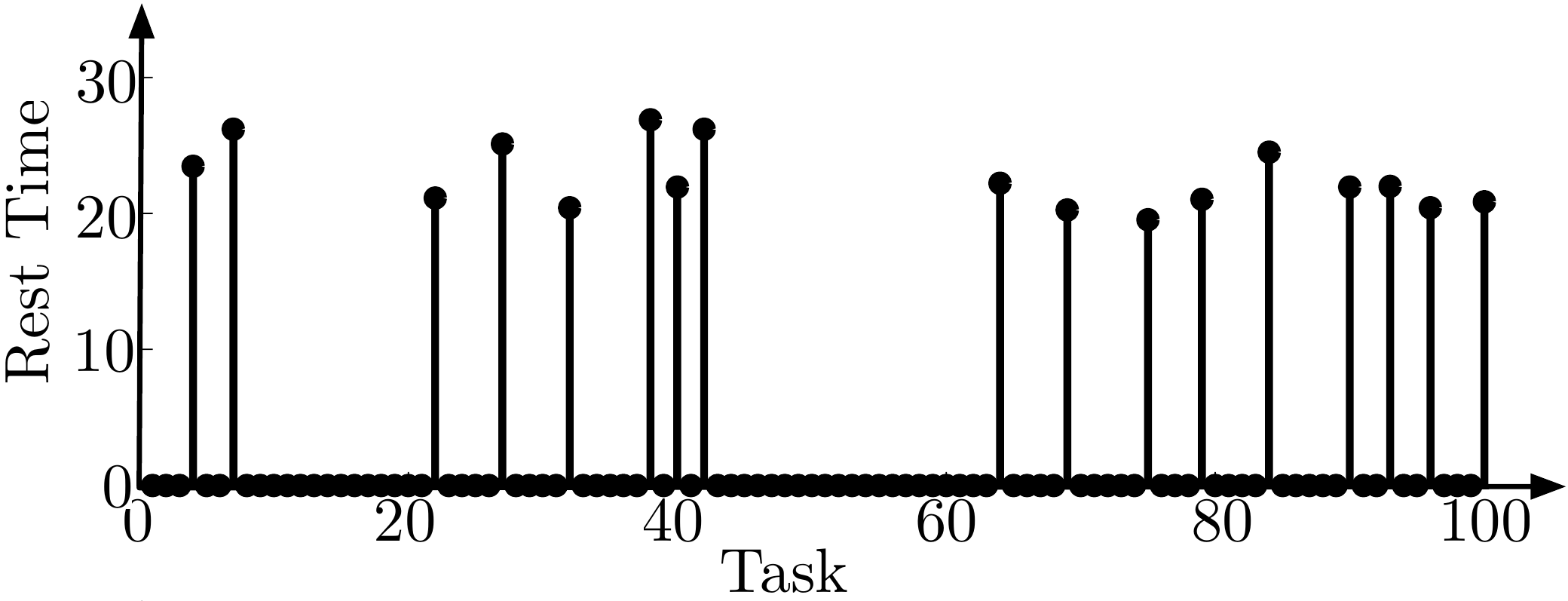} \label{fig:rest} }\\
\subfigure[Queue length]{
\includegraphics[width=0.4\linewidth]{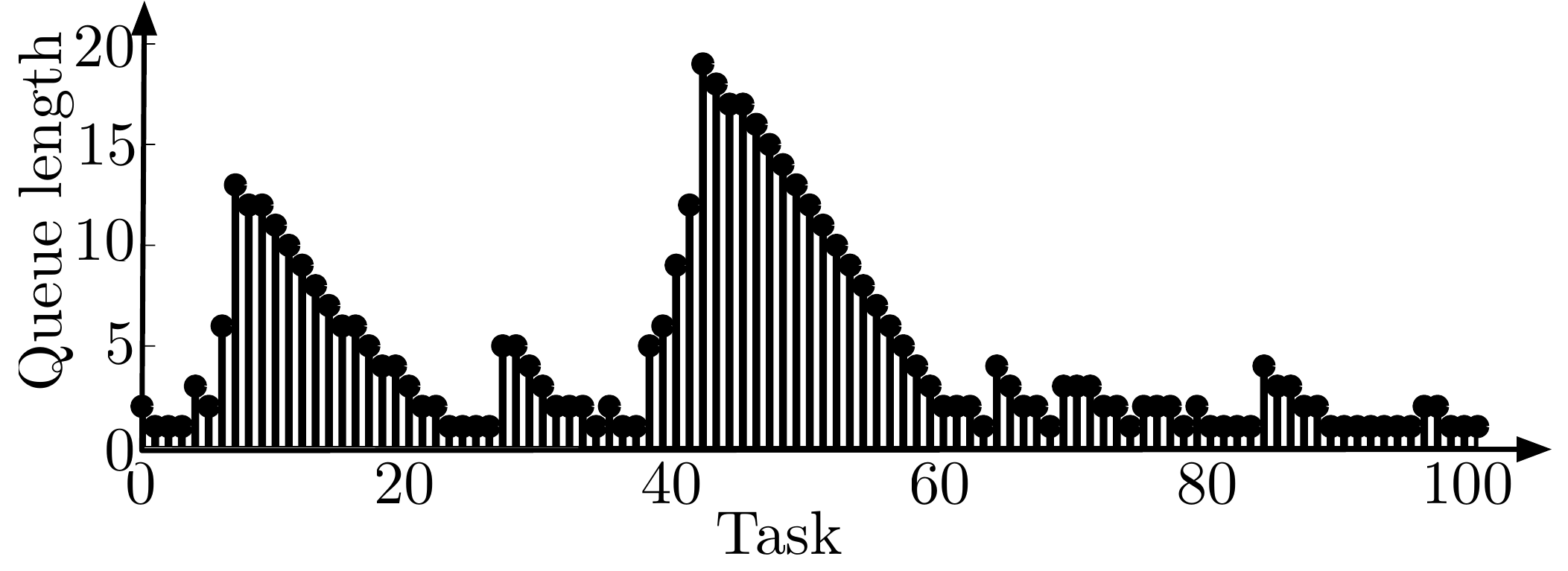}\label{fig:queue-full}} \qquad
\subfigure[Utilization]{
\includegraphics[width=0.4\linewidth]{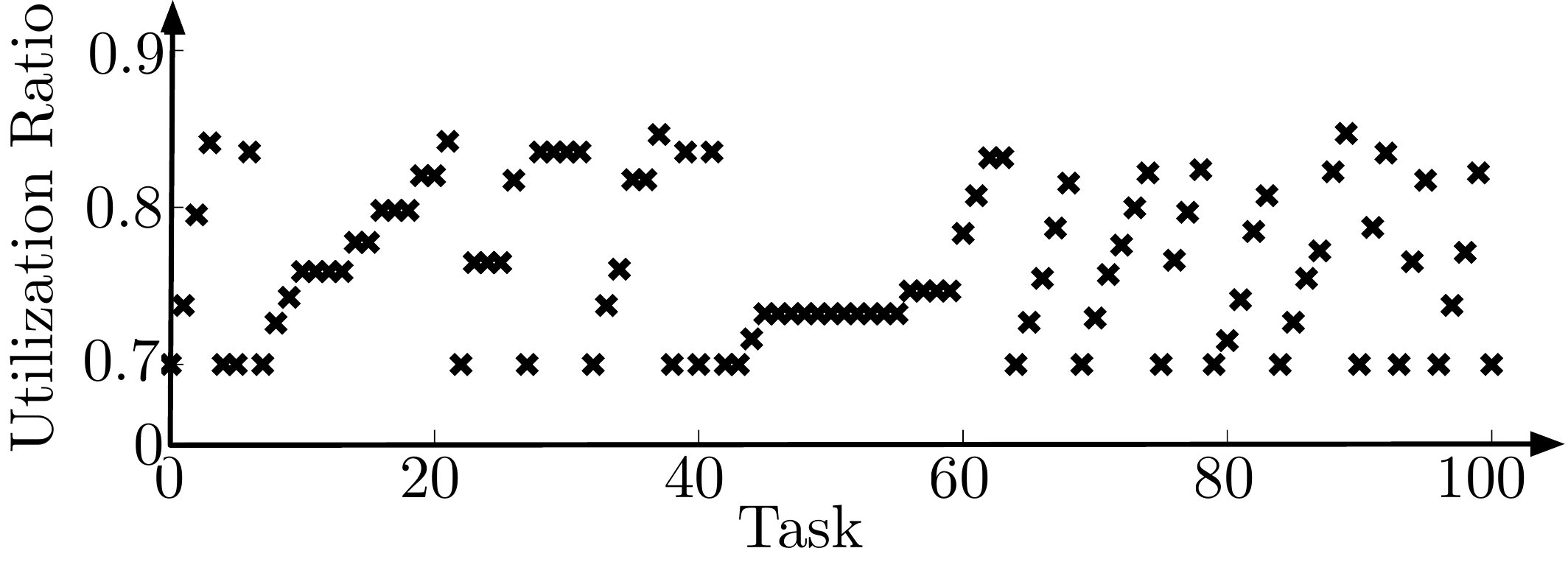}\label{fig:utilization}} \\
\subfigure[Motor Time]{
\includegraphics[width=0.4\linewidth]{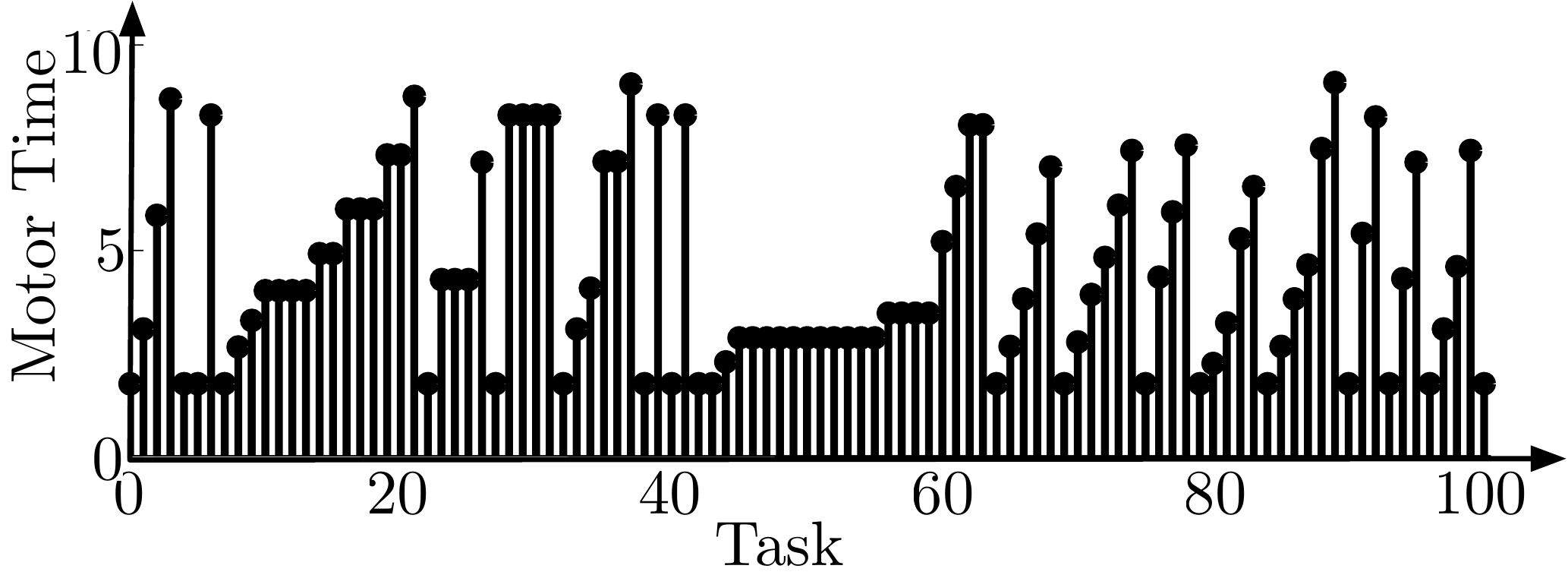}\label{fig:motor}} \qquad
\subfigure[CUSUM statistics]{
\includegraphics[width=0.4\linewidth]{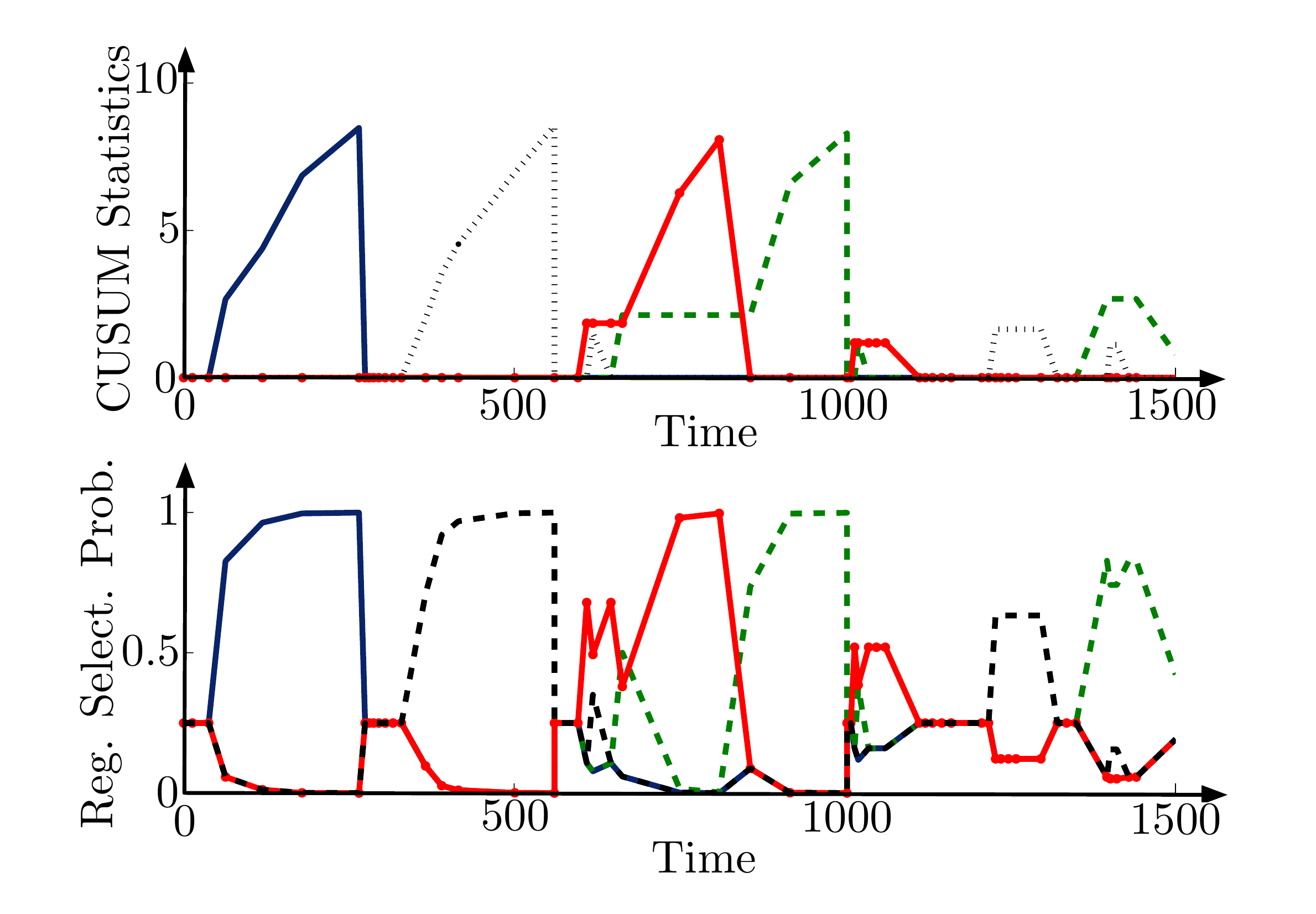}\label{fig:stat-full}}  \\
\subfigure[Region selection probabilities]{
\includegraphics[width=0.4\linewidth]{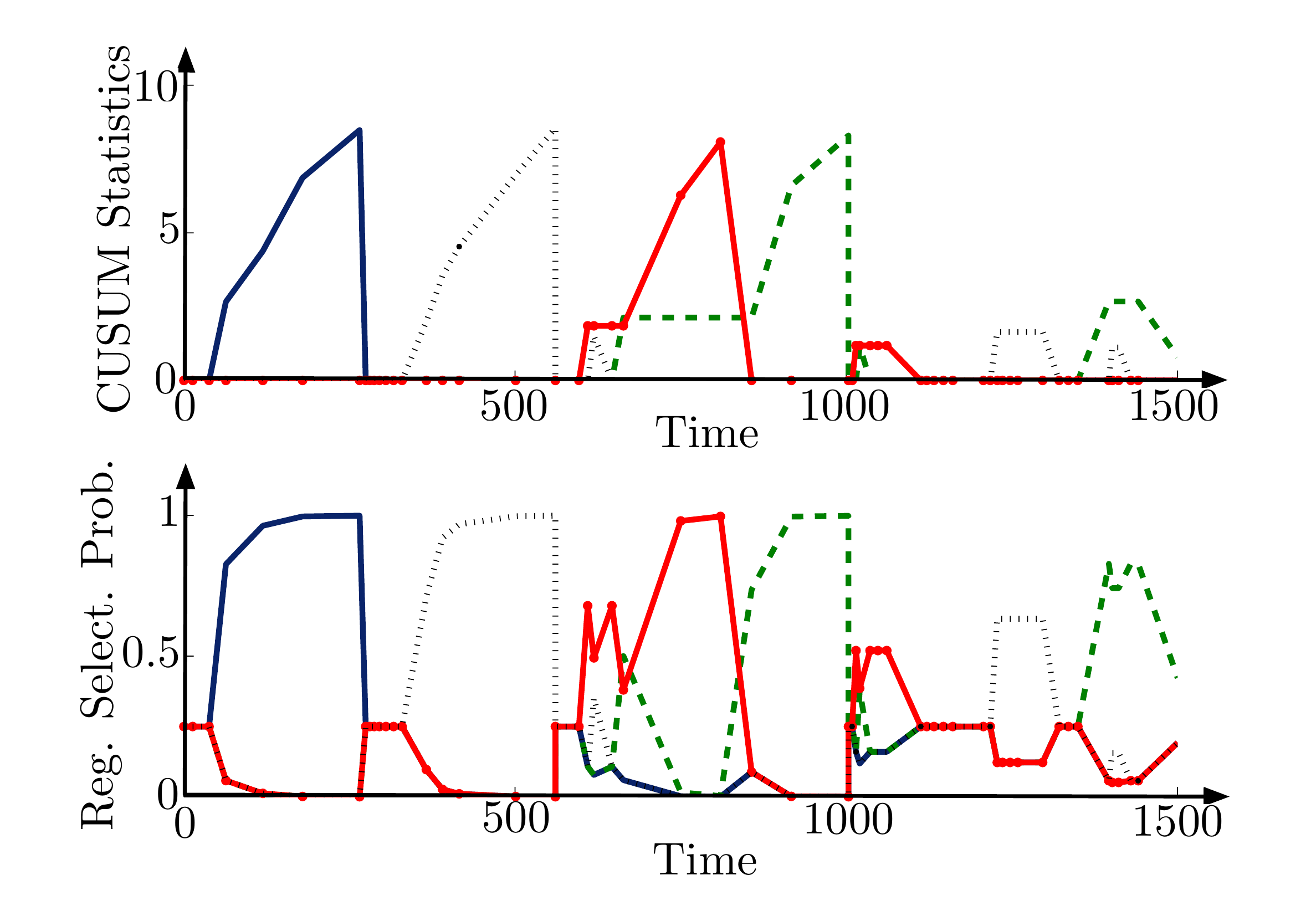}\label{fig:select-full}} \qquad
\subfigure[Retained Belief]{
\includegraphics[width=0.4\linewidth]{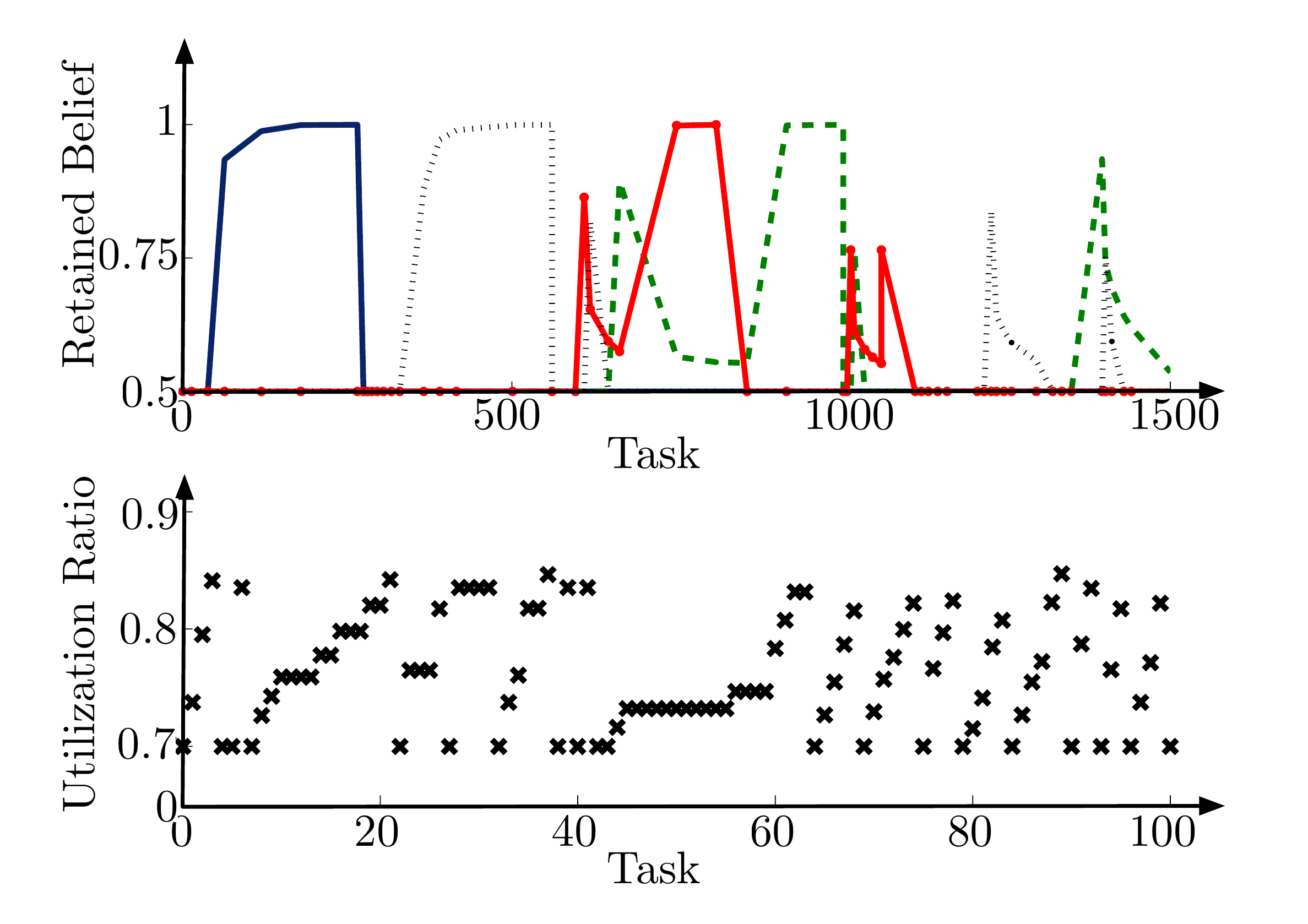}\label{fig:retained-belief}}
\caption{\footnotesize A sample evolution of the CAMS with exogenous factors for a surveillance mission involving four regions. The solid blue lines, dashed green lines, dotted black lines, and solid red lines with dots represent statistics, region selection probabilities, and retained beliefs associated with each region. \label{fig:CAMS-full}}
\end{figure*}

The optimization problem~\eqref{eq:maximize-receding-horizon-real-time} with a horizon length $N=5$
and the modified performance function in~\eqref{eq:certain-performance} is solved before processing each task to determine efficient allocations for the human operator. A sample evolution of the CAMS is shown in Figure~\ref{fig:CAMS-full}. The allocation policy, the queue length, CUSUM statistics, and the region selection policy evolve similarly to the case without exogenous factors and are shown in Figures~\ref{fig:allocation-full},~\ref{fig:queue-full},~\ref{fig:stat-full}, and~\ref{fig:select-full}, respectively. 
The evolution of the operator utilization, associated motor times, and the rest times suggested to bring the operator utilization to an optimal value are shown in Figures~\ref{fig:utilization},~\ref{fig:motor}, and~\ref{fig:rest}, respectively. 
The value of $\subscr{u}{th}$ is set to $0.85$, once the operator utilization crosses $\subscr{u}{th}$, then an appropriate rest time is allocated to the operator such that their utilization drops to $\subscr{u}{opt}=0.7$.
The evolution of the retained belief of the operator is shown in Figure~\ref{fig:retained-belief}. It can be seen that the memory effects take operator's belief about the presence of an anomaly at a region to $0.5$, if that region is not visited for a long time. 

\vspace{-0.1in}
\section{Conclusions}
In this paper, we surveyed models for human decision-making performance and presented a framework for the design of mixed human-robot teams. We demonstrated the framework using the context of a surveillance mission. Along with the fundamental human decision-making performance, we incorporated several other human factors, including, situational awareness, fatigue, and memory retention into human performance model and characterized the accuracy of human decisions in a two-alternative choice task. We utilized these models to design a mixed human-robot team surveillance mission. In particular, we designed simultaneous vehicle routing policies for the vehicles, attention allocation policies for the operator, and anomaly detection algorithms that utilize potentially erroneous decisions by the operator.  

\section{Open Directions}
There are several open directions of research in mixed human-robot teams. Indeed, an obvious open direction is the experimental implementation and testing of the proposed methodology. However, in the following, we focus on the theoretical challenges in the mixed team design.

\paragraph{Exploration-exploitation trade-off}
The system designed relies heavily on the
estimated values of the parameters.  In this setting, real-time evaluation
of parameters becomes critical. One naive strategy for such an evaluation is to introduce
some "control" tasks. For instance, at each stage of the queue, each task may
be a control task with small probability. The response to such control
tasks may be used to update the model parameters in a Bayesian
setting. Moreover, the probability that the current task is a control task
may be adapted depending on the deviation of the operator performance from
the predicted performance.

\paragraph{Multiple anomalies at each region}
In this paper, we assumed fixed drift and diffusion rates for the DDM associated with anomalous as well as non-anomalous regions.This corresponds to a situation in which only one type of anomaly may appear at the region. In general, there may be a set of anomalies that may appear at a region. Multiple anomalies can be modeled using a DDM  in which  drift and diffusion rates are random variables. Different values of these  drift  and diffusion rates correspond to different type of anomalies and different number of anomalies that may be present at the region. As a task is processed, the distributions of drift and diffusion rates can be adapted  according to the operator performance and the GLR algorithm~\cite{MB-IVN:93} can be used as an anomaly detection algorithm.  

\paragraph{Scheduling problems under the free response regime for the operator}
The policies designed in this paper consider the interrogation paradigm of human decision-making. Another important scenario is when the operator makes decision under the free response paradigm. In such a setting, the task scheduling problem becomes important. The processing time associated with each task is a random variable whose distribution depends on several factor including the situational awareness, the fatigue, and  the forgetfulness of the operator.
The scheduling problem involving tasks with stochastic processing time is  a computationally hard problem and the addition of other dynamic elements make it harder. It is of interest to develop approximation algorithms for such scheduling problems.  

\paragraph{Capturing operator-automaton trust}
In this paper, we assumed that the operator trusts the automaton and allocates the duration suggested by the automaton. In general, the operator may not completely trust the automaton and with certain probability may process the task using the free response paradigm. In this setting, the duration allocation to each task is equal to the random processing time associated with the free response paradigm with some probability and is equal to the deterministic time suggested by the decision support system otherwise. It is of interest to determine  optimal policies in this partial trust regime.

\paragraph{Incorporating richer models for human decision-making}
In this paper, we relied on models for human decision-making in two-alternative choice tasks. While the two-alternative choice task is an appropriate model for the human activity in a surveillance mission, other missions may require models for human decision-making in a broader set of tasks, e.g., the multi-armed bandit task~\cite{PR-VS-NEL:13d}. It is of interest to incorporate models for human decision-making in such tasks into the mixed team design. 

\paragraph{Other queuing disciplines and tasks with long reaction times}
In this paper we considered the \emph{first-come-first-serve} processing discipline for the queue of decision-making tasks. It is of interest to study optimal policies for other queueing disciplines, e.g., preemptive priority queues, and the scenario in which tasks in the queue are reordered before the next task is processed. 
In this paper we assumed that the task duration is small and the operator parameters remain constant over the duration of the task. It is of interest to relax this assumption.

\paragraph{Real-time parameter update using observable human actions and states}

An ambitious open direction is to update parameters of the operator model using her observable actions. For instance, the pupil movement of the operator may be monitored and can be used to estimate the parameters for the  cognitive model of the operator. A plausible approach to this end may involve  hidden variable estimation and Bayesian methods to update the parameters of the operator model in real time.

\section*{Acknowledgments}
This work has been supported in part by AFOSR MURI Award FA9550-07-1-0528.

\small
%\bibliographystyle{unsrt}
%\bibliography{alias,FB,Main}

\end{document}